\documentclass[12pt]{article}

\usepackage{amsmath}
\usepackage{mathtools}
\usepackage{amssymb}
\usepackage{here}
\usepackage{mathrsfs}
\usepackage{enumerate}
\usepackage{multirow}
\usepackage{color}
\usepackage{cite}
\usepackage{dsfont}
\usepackage{fullpage}

\newcommand{\pd}{\ensuremath\mathbb{S}_{\succ0}}
\newcommand{\psd}{\ensuremath\mathbb{S}_{\succeq0}}

\everymath{\displaystyle}

\definecolor{lblue}{RGB}{0,110,152}

\definecolor{dred}{RGB}{171,67,53}

\definecolor{dgreen}{RGB}{0,100,0}

\providecommand{\blue}[1]{\color{black}{#1}\color{black}\hspace{0pt}}

\newtheorem{theorem}{Theorem}

\newtheorem{define}[theorem]{Definition}

\newtheorem{hyp}[theorem]{Assumption}
\newtheorem{remark}[theorem]{Remark}

\newenvironment{proof}{{\it Proof :~}}{\hfill$\diamondsuit$\\}

\DeclareMathOperator*{\col}{col}

\DeclareMathOperator{\He}{Sym}
\DeclareMathOperator{\Sym}{Sym}
\DeclareMathOperator*{\diag}{diag}

\DeclareMathOperator{\eps}{\varepsilon}

\DeclareMathOperator*{\argmin}{argmin}
\DeclareMathOperator*{\argmax}{argmax}

\def\P{\mathbb{P}}

\def\E{\mathbb{E}}

\def\d{\mathrm{d}}
\def\dt{\mathrm{d}t}
\def\ds{\mathrm{d}s}

\providecommand\net{$(\Xz,\mathcal{R})\ $}

\def\HJB{Hamilton-Jacobi-Bellman }

\providecommand\X[1]{\boldsymbol{X_{#1}}}

\def\llongrightarrow{\relbar\joinrel\relbar\joinrel\relbar\joinrel\rightarrow}

\providecommand{\rarrow}[1]{\stackrel{#1}{\llongrightarrow}}

\def\Xz{\boldsymbol{X}}

\providecommand\X[1]{\boldsymbol{X_{#1}}}

\def\T{\mathrm{{\scriptstyle T}}}

\begin{document}

\title{Optimal and $H_\infty$ Control of Stochastic Reaction Networks\thanks{This work was supported by the European Research Council under the European Union’s Horizon 2020 research and innovation programme / ERC Grant agreement under Grant 743269 (CyberGenetics)}}
\author{Corentin Briat and Mustafa Khammash\thanks{email: corentin@briat.info, mustafa.khammash@bsse.ethz.ch; url: www.briat.info, https://www.bsse.ethz.ch/ctsb.}\\D-BSSE, ETH-Z\"urich}
\markboth{IEEE TRANSACTIONS ON AUTOMATIC CONTROL, VOL. XX, NO. X, XX XXXX}{Optimal and $H_\infty$ Control of Stochastic Reaction Networks}
\maketitle
\begin{abstract}
Stochastic reaction networks is a powerful class of models for the representation a wide variety of population models including biochemistry. The control of such networks has been recently considered due to their important implications for the control of biological systems. Their optimal control, however, has been relatively few studied until now. The continuous-time finite-horizon optimal control problem is formulated first and explicitly solved in the case of unimolecular reaction networks. The problems of the optimal sampled-data control, the continuous $H_\infty$ control, and the sampled-data $H_\infty$ control of such networks are addressed next. The results in the unimolecular case take the form of nonstandard Riccati differential equations or differential Lyapunov equations coupled with difference Riccati equations, which can all be solved numerically by backward-in-time integration.\\

\end{abstract}

%

\section{Introduction}

\blue{Reaction networks are a very general class of systems that can represent a wide variety of real-world processes including population, ecological, biochemical, epidemiological, and social systems, etc. -- see \cite{Goutsias:13} for a survey. Deterministic reaction network are usually represented in terms differential equations and their extensions such as delay-differential equations or partial differential equations, when delays or diffusion effects need to be taken into account. Tools for analyzing and controlling them are now well established even though some more specialized tools, tailored to the framework of reaction networks, have also been derived for analyzing their stability using, for instance, chemical reaction network theory \cite{Horn:72,Feinberg:72} or Lyapunov methods \cite{Blanchini:14,AlRadhawi:16,AlRadhawi:20}. The control of such networks is now a growing interest in the control community  \cite{Briat:15e,Cuba:17,Xiao:18,Briat:19:Logistic,Cuba:19cdc,Agrawal:19cdc,DelVecchio:15,DelVecchio:16,DelVecchio:17,Khammash:19,Olsman:19,Olsman:19b}.}


\blue{In a biological context, deterministic models are only accurate when the populations are both large in size and homogenously distributed so that the interactions between the different network components can be assumed to evolve in a continuous fashion and the network participants can be characterized in terms of a continuous variables; such as concentrations in chemistry. However, when the network participants are in low-copy numbers, the assumption of homogeneity may break down, the actual count of the network participants needs now to be kept track of, and the randomness in the reactions can not be neglected anymore. In such a case, it is necessary to consider stochastic reaction networks which have been shown to be of crucial importance in biology due to their ability to capture both intra- and inter-cellular variability \cite{Elowitz:02}. Such networks can then be represented as continuous time Markov jump processes \cite{Anderson:15} and their approximations. The analysis of stochastic reaction networks is still an open problem for which partial solutions have been obtained and formulated in terms of both theoretical and computational tools \cite{Briat:13i,Engblom:12,Rathinam:13}. Theoretical results regarding their control have also been obtained in \cite{Briat:12c,Briat:13h,Briat:19:Opto} in the moment equations framework and in \cite{Briat:15e,Briat:19:DelayedRN,Briat:18:Interface,Gupta:19,Briat:20:Structural} in the framework of chemical reaction network controllers, that is, controllers implemented as stochastic reaction networks themselves. Practical implementations of such control strategies in a biological setting have been reported in  \cite{Milias:11,Uhlendorf:12,Milias:16,Fiore:16,Lugagne:17,Shannon:20}, and in \cite{Fiore:16,Aoki:19,Shannon:20} in the contexts of \emph{in-silico} controllers  and \emph{in-vivo} controllers, respectively.}

\blue{Until now, all the designed controllers were chosen such that some simple control objectives were satisfied, such as set-point regulation, constant disturbance rejection, robustness to uncertainties, robustness to noise, or even noise exploitation. In this regard, the controller structures able to ensure (robust) perfect adaptations were restricted to Proportional-Integral-Derivative-like controllers and their approximations. While this is a perfectly viable solution in the Single-Input-Single-Output (SISO) case, the Multiple-Input-Multiple-Output (MIMO) case makes the tuning of such controllers a much harder task. A natural step forward is, therefore, the design of controllers that would minimize some cost in order to solve optimal MIMO control problems -- as historically done in the field of control with the introduction of the Linear Quadratic (LQ) and the Linear Quadratic Gaussian (LQG) control, to cite a few. The two main analytical approaches for solving optimal control problems are Pontryagin's Maximum Principle and Dynamic Programming. The former is essentially based on the calculus of variations, which is an extension of the usual Lagrange multiplier approach for standard optimization \cite{Boyd:04,Luenberger:69,Ross:09}. It can be seen as a trajectory-based (or temporal-based) approach and is convenient to use whenever state and control constraints are involved. The latter, on the other hand, is based on Bellman's Principle of Optimality \cite{Bellman:54}, and unlike the Maximum Principle, it is a spatial approach based on the use of a so-called \emph{value function} which needs to satisfy a certain optimality condition at every point in the state-space \cite{Bertsekas:17,Bertsekas:12}. Both approaches have benefits and drawbacks, and their connections are now well understood.}

\blue{The objective of this paper is to develop an optimal control theory for stochastic reaction networks. This problem, despite its importance, as been relatively few studied until now. General results on the optimal control of general Markov processes have been obtained for instance in \cite{Stockbridge:90,Kurtz:98}. The optimal control of stochastic reaction networks is considered in \cite{Zhang:18} under the assumption that the copy number of the involved species is large, which justifies the use of limiting arguments, and that the set of possible actions is finite. Some approximate online algorithms are proposed there based on state-space truncation and on the empirical evaluation of the cost to optimize. The optimal control of infinite-dimensional piecewise deterministic Markov processes is considered in \cite{Renault:17} with some application to the control of neuronal dynamics via optogenetics. In this setup, the process switches from one dynamics to another following a finite Markov chain whose rates are control input. An alternative approach for the control of stochastic reaction networks is to first approximate the Markov jump process as a diffusion process \cite{VanKampen:07}, for instance, on which standard tools such as Pontryagin's Maximum Principle and Dynamic Programming can be applied.  However, it is known that those approximate dynamics may dramatically differ from that of the original jump Markov process and, in this regard, this does not seem to a valid option. }


\blue{We propose in this paper to keep the original Markov jump process formulation and to address the optimal control problem using Dynamic Programming directly. Even though it is known that such an approach suffers from the so-called curse-of-dimensionality, it seems interesting, at least, to develop it for theoretical purposes. Similar approximate solving methods to those that can be used for standard optimal control problems can be considered in the present context. For this reason, they will only be briefly mentioned here without entering into details. The optimal control of Markov chains and Markov processes has also been addressed but essentially only in the finite-state space setting \cite{Bertsekas:17,Bertsekas:12,Renault:17}; an approach that has also been adapted to the infinite state-space case through either its discretization or its truncation \cite{Zhang:18}, as often done in the Reinforcement Learning and Control of Markov Decision Processes \cite{Bertsekas:19}. The approach we consider here keeps the structure of the problem intact and deals with the actual system and its state-space structure directly.}

\blue{We first address the continuous-time finite-horizon optimal control of stochastic reaction networks and provide the characterization of the optimal continuous-time control law in terms of a Hamilton-Jacobi-Bellman (HJB) equation taking, in this case, the form of a differential-difference equation. This has to be contrasted with the usual partial differential equations obtained in optimal control and the fact that viscosity solutions do not play an as important role as in previous works \cite{Fleming:06}. In fact, it should only be necessary to consider mild-solutions instead. In the special case of unimolecular networks, the HJB equation exactly reduces to a nonstandard Riccati differential equation which does not seem to have been previously obtained in the literature; see e.g. \cite{AbouKandil:03}. However, its structure is sufficiently similar to existing ones so that some of the existing analysis tools can be considered; e.g. for the establishing the existence of solutions. Unfortunately, it is not immediate to see any connection with a "Hamiltonian matrix" as in previous works on the optimal control of deterministic or diffusion processes.}

\blue{We provide then a solution to the optimal sampled-data control of stochastic reaction networks. To this aim, we first reformulate the overall system as a stochastic hybrid system \cite{Teel:14} which contains both random and deterministic jumps. While the former correspond to state updates from the randomly firing reactions, the latter corresponds to the control input updates, which are assumed to occur periodically. The advantage of this reformulation lies in the fact that the structure of the system remains unchanged unlike in other methods such as in discretization-based ones. Note also that the overall procedure to exactly discretize a stochastic reaction network with infinite state-space is unclear and, even if possible, it would result in a dramatic reformulation of the model of the system which would not be exploitable for control purposes. A final drawback of discretization-based methods is the loss of the so-called \textit{inter-sample behavior} which needs is involved in the cost expression. This issue notably motivated the introduction of lifting methods \cite{Yamamoto:90,Yamamoto:94} and other subsequent ones. Based on the hybrid reformulation, we characterize the optimal control law in terms of hybrid HJB equations. In the unimolecular case, those equations reduce to a Lyapunov differential equation coupled with a discrete-time Riccati difference equation.}

\blue{We next remove the assumption that some of the uncontrolled inputs are known and consider them as unknown exogenous disturbances. To account for this, we consider an $H_\infty$ performance criterion, introduced in the control literature in \cite{Zames:81} and refined in \cite{Doyle:89a,Basar:95} in the context of optimal control and game theory. Such a framework allows to quantify the impact of those exogenous inputs on the so-called controlled/performance outputs in an $L_2$-sense and from which we develop an approach to compute the optimal control law that minimizes the worst-case $L_2$-gain of the transfer from exogenous inputs to the controlled outputs. This is achieved using a game-theoretic approach which exactly characterizes both the optimal control law and the worst-case exogenous inputs in terms of the solution of a Hamilton-Jacobi-Isaacs (HJI) equation \cite{Isaacs:65}. In the case of stochastic reaction networks, this equation takes the form of a nonlinear differential-difference equation which reduces, in the unimolecular case, to a non-standard Riccati differential equation which does not seem to have been reported nor studied elsewhere before.}

\blue{As in the continuous-time optimal control problem, we also extend those results to the sampled-data case where we use, again, a hybrid formulation of the system which allows to keep all the necessary information about the system between the sampling instants. The discretization of the process would result in the loss of both the structure of the system and its inter-sample behavior while the use of lifting methods would require the consideration of integral input and output operators \cite{Chen:95}, which would dramatically complicate the overall approach. The hybrid approach circumvents all those issues in an elegant way and allows to characterize the optimal sampled-data control law and the worst-case exogenous disturbances in terms of the solution of hybrid HJI equations which reduce, in the unimolecular case, to a Riccati differential equation coupled with a discrete-time Riccati difference equation.}\\

\blue{\noindent\textbf{Notations.} The cone of real positive (semi)definite matrices of dimension $n$ is denoted by ($\mathbb{S}^n_{\succeq0}$) $\mathbb{S}^n_{\succ0}$. For a square real matrix $A$, we denote $\He[A]:=A+A^{\T}$. For a differentiable function $V(t,x)$, we denote by $V_t(t,x)$ its partial derivative with respect to $t$. For some scalars $x_1,\ldots,x_n$, the vector consisting of stacking those entries vertically with $x_1$ on top is denoted by $\textstyle\col_{i=1}^n(x_i)$ or $\col(x_1,\ldots,x_n)$. We denote by $e_{d+1}$ the vector of zeros of dimension $n$ ($n\ge d+1$) except at the position $d+1$ where the entry is one.}\\

\blue{\noindent\textbf{Outline.} Preliminaries on stochastic reaction networks are recalled in Section \ref{sec:preliminaries}. The optimal continuous-time control of stochastic reaction networks and its sampled-data counterpart are considered in Section \ref{sec:CT:LQ} and Section \ref{sec:SD:LQ}, respectively. Similarly, the optimal continuous-time $H_\infty$ control of stochastic reaction networks and its sampled-data counterpart are addressed in Section \ref{sec:CT:LQ} and Section \ref{sec:SD:Hinf}, respectively.}

\section{Preliminaries}\label{sec:preliminaries}

\subsection{Stochastic Reaction Networks}\label{sec:SRN}

\blue{A reaction network $(\Xz,\mathcal{R})$ consists of a set of $d$ molecular species $\Xz=\{\X{1},\ldots,\X{d}\}$ that interacts through $K$ reaction channels $\mathcal{R}=\{\mathcal{R}_1,\ldots,\mathcal{R}_K\}$ denoted as
\begin{equation}
 \mathcal{R}_k:\ \sum_{i=1}^d\zeta_{k,i}^l\X{i}\rarrow{\rho_k(t)}  \sum_{i=1}^d\zeta_{k,i}^r\X{i},\ k=1,\ldots,K
\end{equation}
where $\rho_k:\mathbb{R}_{\ge0}\mapsto\mathbb{R}_{\ge0}$ is the (time-varying) reaction rate parameter and $\zeta_{k,i}^\ell,\zeta_{k,i}^r\in\mathbb{Z}^d_{\ge0}$ are the left and right stoichiometric vectors. The stoichiometric vector of reaction  $\mathcal{R}_k$ is given by $\zeta_k:=\zeta_k^r-\zeta_k^\ell\in\mathbb{Z}^d$ where $\zeta_k^r=\col(\zeta_{k,1}^r,\ldots,\zeta_{k,d}^r)$ and $\zeta_k^l=\col(\zeta_{k,1}^l,\ldots,\zeta_{k,d}^l)$. That is, when the reaction $\mathcal{R}_k$ fires, the state jumps from $x$ to $x+\zeta_k$. The stoichiometry matrix $S\in\mathbb{Z}^{d\times K}$ is defined as $S:=\begin{bmatrix}
  \zeta_1\ldots\zeta_K
\end{bmatrix}$. When  the kinetics is  mass-action, the propensity function of reaction $\mathcal{R}_k$ is given by $\lambda_k(t,x)=\rho_k(t)P_k(t,x)$ where
\begin{equation}\label{eq:massactionpoly}
  P_k(t,x)=\prod_{i=1}^d\begin{pmatrix}
    x_i\\
    \zeta_{k,i}
  \end{pmatrix}
\end{equation}
and is such that  $\lambda_k(t,x)=0$ if $x\in\mathbb{Z}_{\ge0}^d$ and $x+\zeta_k\notin\mathbb{Z}_{\ge0}^d$. Under the well-mixed assumption, this network can be described by a continuous-time Markov process $(X_1(t),\ldots,X_d(t))_{t\ge0}$ with state-space $\mathbb{Z}_{\ge0}^d$; see e.g. \cite{Anderson:11}.

For any locally bounded function $V:\mathbb{R}\times\mathbb{Z}^d_{\ge0}\mapsto\mathbb{R}$, the generator $\mathbb{A}$ of the Markov process associated with the reaction network $(\Xz,\mathcal{R})$ is defined as
\begin{equation}
  (\mathbb{A}V)(t,x)=V_t(t,x)+\sum_{i=1}^K\lambda_i(t,x)\Delta_iV(t,x)
\end{equation}
where $\Delta_iV(t,x):=V(t,x+\zeta_i)-V(t,x)$.}


\subsection{Kolmogorov's equations}\label{sec:Kolmogorov}

\blue{\noindent\textbf{Forward Kolmogorov Equation.} Let $(X_1(t),\ldots,X_d(t))_{t\ge0}$ be the CTMC representing the stochastic reaction network $(\Xz,\mathcal{R})$ with some initial state $X(0)\in\mathbb{Z}_{\ge0}^d$. For any state value $x\in\mathbb{Z}_{\ge0}^d$, we can define the probability to be in that state at time $t$ as
\begin{equation}
  p(t,x) := \P(X(t)=x|X(0)).
\end{equation}
It is known that the time evolution of this probability distribution is described by the so-called \textit{Forward Kolmogorov Equation} or \textit{Chemical Master Equation} (CME) \cite{Anderson:11} given by
\begin{equation}\label{eq:CME}
  p_t(t,x)=\sum_{k=1}^K\lambda_k(t,x-\zeta_k)p(t,x-\zeta_k)-\sum_{k=1}^K\lambda_k(t,x)p(t,x),\ x\in\mathbb{Z}_{\ge0}^d.
\end{equation}
The CME is a linear differential equation of dimension equal to the number of elements in the state-space. In the present case, this number is infinite and of the cardinality of the natural numbers. The CME is, in general, difficult to solve analytically and this is the reason why only few analytical results exist.\\



\noindent\textbf{Backward Kolmogorov Equation.} Let $f:\mathbb{Z}^d_{\ge0}\mapsto\mathbb{R}$ be a locally bounded function. Then, the backward time evolution of the function $V(t,x) := \E[f(X(T))|X(t)=x)$, with $\E[f(X(T))]$ finite, according to the CTMC $(X_1(t),\ldots,X_d(t))_{t\ge0}$ representing the stochastic reaction network $(\Xz,\mathcal{R})$ over the time interval $[0,T]$ is governed by the \textit{Backward Kolmogorov equation}
\begin{equation}
  -V_t(t,x)=\sum_{k=1}^K\lambda_k(t,x)(V(t+\zeta_k)-V(t,x))
\end{equation}
together with the terminal condition $V(T,x)=f(x)$. This equation will turn out to play a crucial role in the optimal control of stochastic reaction networks. 
}

\subsection{Moment equations}\label{sec:moments}

\blue{Based on the CME, dynamical expressions for the first- and second-order moments may be easily derived and are given by
\begin{equation}\label{eq:moments_g}
\hspace{-3mm}\begin{array}{lcl}
    \dfrac{\d \E[X(t)]}{\dt}&=&S\E[\lambda(t,X(t))],\\
    \dfrac{\d \E[X(t)X(t)^{\T}]}{\dt}&=&\He[S\E[\lambda(t,X(t))X(t)^{\T}]]+S\diag\{\E[\lambda(t,X(t))]\}S^{\T}.
\end{array}
\end{equation}
Those equations cannot be directly solved in the general case due to some moment closure issues (see e.g. \cite{Singh:11}) but become exactly solvable in the unimolecular case where the propensity functions take the form $\lambda(t,x)=W_x(t)x+W_0(t)$, for some suitable functions $W_x:\mathbb{R}_{\ge0}\mapsto\mathbb{R}_{\ge0}^{K\times d}$, $\textstyle W_x(t)=:\col_{i=1}^K(w_{x,k}^{\T})$, $W_0:\mathbb{R}_{\ge0}\mapsto\mathbb{R}_{\ge0}^K$, $\textstyle W_0(t)=:\col_{i=1}^K(w_{0,k}^{\T})$. In such a case, the moment equations can be rewritten as
\begin{equation}\label{eq:moments}
\begin{array}{rcl}
    \dfrac{\d \E[X(t)]}{\dt}&=&A(t)\E[X(t)]+d(t),\\
     \dfrac{\d U(t)}{\dt}&=&\He\left(A(t)U(t)+d(t)\E[X(t)^{\T}]\right)+S\diag\{W_x(t)\E[X(t)]+W_0(t)\}S^{\T}
\end{array}
\end{equation}
where $U(t):=\E[X(t)X(t)^{\T}]$, $\textstyle A(t)=SW_x(t)=\sum_{i=1}^n\zeta_kw_{x,k}^{\T}$, $\textstyle d(t)=SW_0(t)=\sum_{i=1}^n\zeta_kw_{0,k}^{\T}$. Introducing $\bar{X}(t)=\col(X(t),1)$ and $\bar U(t):=\E[\bar X(t)\bar X(t)^{\T}]$, we get that
\begin{equation}\label{eq:moments2}
\begin{array}{rcl}
    \dfrac{\d \E[\bar X(t)]}{\dt}&=&\bar A(t)\E[\bar X(t)],\\
     \dfrac{\d \bar U(t)}{\dt}&=&\bar A(t)\bar U(t)+\bar U(t)\bar A(t)^{\T}+\dfrac{1}{2}\sum_{i=1}^K\left(\bar{A}_k(t)\bar U(t)e_{d+1}\zeta_k^{\T}+\zeta_ke_{d+1}^{\T}\bar U(t)\bar{A}_k(t)^{\T}\right)
\end{array}
\end{equation}
where
\begin{equation}
  \bar A(t)=\begin{bmatrix}
    A(t) & d(t)\\
    0 & 0
  \end{bmatrix},\ \bar A_k(t)=\begin{bmatrix}
    \zeta_kw_{x,k}(t) & \zeta_kw_{0,k}(t)\\
    0 & 0
  \end{bmatrix}.
\end{equation}
It is interesting to note that the dynamics of the mean molecular species is described by a Linear Time-Varying (LTV) positive system, meaning that state is confined in the nonnegative cone \cite{Farina:00}. This can be verified from the fact that the matrix $A(\cdot)$ is Metzler (i.e. all the off-diagonal entries are nonnegative and $d(\cdot)$ is nonnegative. Similarly, the second-order moments are confined in the positive semidefinite cone, which is a well-known property.}

\section{Continuous-Time Optimal Control of Stochastic Reaction Networks}\label{sec:CT:LQ}

\blue{The objective of this section is the development of optimal control results for stochastic reaction networks. A general result is given first, and adapted to the case of stochastic reaction networks with mass-action kinetics. The latter result is then specialized to the case of unimolecular reaction networks having affine propensity functions. This case is analogous to that of the optimal control of finite-dimensional linear time-varying systems and, as such, it will be shown that the solution of the optimal control problem can be formulated in terms of the solution of a differential Riccati equation.}

\subsection{General stochastic reaction networks}

The class of reaction networks considered in Section \ref{sec:preliminaries} are autonomous networks in the sense that those networks do not possess any input channels that can be actuated. In this regard, we need to extend the class of networks of Section \ref{sec:preliminaries} to intercorporate some inputs channels acting at the level of propensity functions; e.g. at the level of the reaction rates. This choice is motivated by cybergenetical or optogenetical applications where networks have been controlled by acting on their propensity functions using an external stimulus such as light \cite{Khammash:11,Lugagne:17,Rullan:17,Renault:17} or chemical inducers \cite{Uhlendorf:12}. In this regard, it seems natural to make the following assumption on the propensity functions:
\begin{hyp}\label{hyp:CTLQ}
  The propensity functions of the network \net are locally bounded functions of the form $\lambda_k(t,x,u)$, $k=1,\ldots,K$, where $u\in\mathbb{R}^m_{\ge0}$ is the control input vector, and, as before, $x\in\mathbb{Z}^d_{\ge0}$ is the state of the network.
\end{hyp}


In order to control the system in the desired way, we introduce the following cost
\begin{equation}\label{eq:costCTgeneral}
\begin{array}{rcl}
  J_T(u)&:=&\E\left[\displaystyle\int_0^Tq(s,X(s),u(s))\ds+q_T(X(T))\right]
\end{array}
\end{equation}
where $q:[0,T]\times \mathbb{Z}_{\ge0}^d\times \mathbb{R}_{\ge0}^m\mapsto\mathbb{R}_{\ge0}$ is the running cost and $q_T:\mathbb{Z}_{\ge0}^d\mapsto\mathbb{R}_{\ge0}$ is the terminal cost. The control input is then chosen such that this cost is minimized
\begin{equation}
  u^*=\argmin_{u\in\mathcal{U}_{0,T}} J_T(u)
\end{equation}
where we have assumed that such a minimum exists and where $\mathcal{U}_{s,t}$ is the set of admissible control inputs defined as
\begin{equation}\label{eq:Ust}
  \mathcal{U}_{s,t}:=\left\{u\in L_2([s,t]\times\mathbb{Z}_{\ge0}^d,\mathbb{R}_{\ge0}^m):\ u\ \textnormal{mesurable}\right\},\ 0\le s\le t\le T.
\end{equation}

In order to solve the above optimization problem, we consider the Dynamic Programming approach \cite{Bellman:57,Bertsekas:17,Bertsekas:12} which relies on the concept of \emph{value function} defined as
\begin{equation}
  V(t,x) := \min_{u\in\mathcal{U}_{t,T}} \E\left[\left.\int_t^Tq(s,X(s),u(s))\ds\right|X(t)=x\right],\ V(T,x)=q_T(x)
\end{equation}
for $x\in\mathbb{Z}_{ge0}^d$ and which leads to the following result:
%
\begin{theorem}
Let us consider a general reaction network \net satisfying Assumption \ref{hyp:CTLQ}. Then, the optimal control law that minimizes the cost \eqref{eq:costCTgeneral} is given by
\begin{equation}
  u^*(t,x)=\argmin_u\left\{\sum_{i=1}^K\lambda_i(t,x,u)\Delta_iV(t,x)+q(t,x,u)\right\}
\end{equation}
  where the function $V(t,x)$ is the solution to the Hamilton-Jacobi-Bellman equation given by
\begin{equation}\label{eq:HBJCT}
  V_t(t,x)+\min_u\left\{\sum_{i=1}^K\lambda_i(t,x,u)\Delta_iV(t,x)+q(t,x,u)\right\}=0,\ V(T,x)=q_T(x)
\end{equation}
for $(x,t)\in\mathbb{Z}_{\ge0}^d\times[0,T]$. Moreover, in such a case, the optimal cost $J_T(u^*)$ coincides with $V(0,X(0))$.
\end{theorem}

\blue{It seems important to stress here that the above equation is vastly different from the standard Hamilton-Jacobi-Bellman equation for the optimal control of continuous-time systems governed by (stochastic) differential equations. In those cases, the value function is described in terms of first-order or second-order elliptic partial differential equations. In the present case, we simply have a differential equation in time whereas the partial derivative operator with respect to the state is replaced by difference operators due to the presence of jump dynamics. In this regard, the solutions are expected to be quite different from those in the standard setup of the control of systems governed by differential equations and stochastic differential equations. In this regard, the concept of viscosity solutions does not play an as essential role as in most of existing optimal control problems.}

\blue{We now focus on the optimal control of stochastic reaction networks with mass-action kinetics as defined in \eqref{eq:massactionpoly} where each one of the reaction rates is now a potential control input. This leads to the consideration of the following assumption:
\begin{hyp}\label{hyp:CTLQ2}
  The propensity functions of the network \net are of the form
  \begin{equation}\label{eq:propensities_bimolecular}
  \lambda_i(t,x,u)=f_i(t,x)+g_{i}(t,x)u,\ i=1,\ldots,K
\end{equation}
where the $f_i(t,x)$'s are suitable polynomials and the $g_{i}(t,x)$'s are suitable row-vector polynomials. Moreover, the function $q(t,x,u)$ in the cost \eqref{eq:costCTgeneral} takes the form
\begin{equation}\label{eq:kdslkd;skdkaldkal}
    q(t,x,u)=\tilde q(t,x)+u^{\T}R(t)u
  \end{equation}
  where $\tilde q(t,x)$ is some nonnegative polynomial and $R(\cdot)\succ0$.
\end{hyp}}

\noindent \blue{We then have the following result:
\begin{theorem}
  Let us consider a reaction network $(\X{},\mathcal{R})$ that satisfies Assumption \ref{hyp:CTLQ2}. Then, the optimal control law minimizing the cost \eqref{eq:costCTgeneral} with \eqref{eq:kdslkd;skdkaldkal} is given by
\begin{equation}\label{eq:optimal_input_CT_bi}
  u^*(t,x)=-\dfrac{1}{2}R(t)^{-1}g(t,x)^{\T}\Delta V(t,x),
\end{equation}
where the value function $V(t,x)$ satisfying $V(T,x)=q_T(x)$, solves the differential-difference equation over $t\in[0.T]$
\begin{equation}\label{eq:HJB_CT_bi}
  V_t(t,x)+\sum_{i=1}^Kf_i(t,x)\Delta_iV(t,x)  -\dfrac{1}{4}\Delta V(t,x)^{\T}g(t,x)R(t)^{-1}g(t,x)^{\T}\Delta V(t,x)+\tilde q(t,x)=0,
\end{equation}
where $\textstyle g(t,x)=\col_{i=1}^K(g_i(t,x))$, and $\textstyle \Delta V(t,x) = \col_{i=1}^K(\Delta_iV(t,x))$.
\end{theorem}
\begin{proof}
  The Hamilton-Jacobi-Bellman equation is given in this case by
\begin{equation}\label{eq:kdlsk;dlksakdd;laksd;l}
  V_t(t,x)+\sum_{i=1}^Kf_i(t,x)\Delta_iV(t,x)+\min_u\left\{\sum_{i=1}^K\Delta_iV(t,x)g_i(t,x)u+q(t,x,u)\right\}=0
\end{equation}
Since  $q(t,x,u)$ is strictly convex in $u$ and $\lambda(t,x,u)$ is affine in $u$, then there exists a unique global minimum to \eqref{eq:HBJCT}. Computing the derivative of the expression in brackets in \eqref{eq:kdlsk;dlksakdd;laksd;l} with respect to $u$ yields the optimal control input given in \eqref{eq:optimal_input_CT_bi} and substituting this expression in  \eqref{eq:HBJCT} yields \eqref{eq:HJB_CT_bi}.
\end{proof}}

\blue{\begin{remark}
  The differential-difference equation \eqref{eq:HJB_CT_bi} is difficult to solve in the general setting. It is, however, possible to find an approximate solution by turning this the equation into an inequality as
\begin{equation}\label{eq:HJB_CT_bi_relax}
  V_t(t,x)+\sum_{i=1}^Kf_i(t,x)\Delta_iV(t,x)  -\dfrac{1}{4}\Delta V(t,x)^{\T}g(t,x)R(t)^{-1}g(t,x)^{\T}\Delta V(t,x)+\tilde q(t,x)\ge0
\end{equation}
together with $V(T,x)\ge q_T(x)$ and by using optimization methods, such as SOS programming \cite{Parrilo:00}, under the assumption that the propensity functions and the approximate value function are polynomial. Note that, In such a case, $V(0,X(0))$ will only be an upper bound on the optimal cost.
\end{remark}}

\blue{\subsection{Unimolecular stochastic reaction networks}

In the case of unimolecular stochastic reaction networks with mass action kinetics, the propensity functions reduce to
\begin{equation}\label{eq:propensities:unimolecular:inputs}
  \lambda_i(t,x,u)= \bar w_i(t)^{\T}\bar x+\bar b_i(t)^{\T}u
\end{equation}
where $\bar x:=\col(x,1)$, where $\bar w_i(t) = \col(w_{i,x}(t),w_{i,0}(t))$, $w_{i,x}\in\mathbb{R}_{\ge0}^{d}$, $\bar w_{i,0}(t)\in\mathbb{R}_{\ge0}$, and $\bar b_i\in\mathbb{R}_{\ge0}^m$, $i=1,\ldots,K$. As already mentioned in Section \ref{sec:moments}, an important property of unimolecular reaction networks is that their moment equation is closed. In the present case, one can write the following expression for the first-order moment equation
\begin{equation}\label{eq:moment:uni:CT}
  \dfrac{\d}{\dt}\E[X(t)]=A(t)\E[X(t)]+B(t)\E[u(t)]+d(t)
\end{equation}
where
\begin{equation}
  A(t)=\sum_{k=1}^K\zeta_k \bar w_{k,x}(t)^{\T},\ B(t)=\sum_{k=1}^K\zeta_k \bar b_k(t)^{\T},\ \textnormal{and }  d(t)=\sum_{k=1}^K\zeta_k \bar w_{k,0}(t).
\end{equation}


As unimolecular reaction networks are analogues of linear systems, the following Linear Quadratic (LQ) cost is considered
  \begin{equation}\label{eq:LQcost}
\begin{array}{rcl}
  J_T(u)&:=&\E\left[\displaystyle\int_0^T\left(\bar{X}(s)^{\T}Q(s)\bar{X}(s)+u(s)^{\T}R(s)u(s)\right)\ds+\bar X(T)^{\T}Q_T\bar X(T)\right]
\end{array}
\end{equation}
where $Q:[0,T]\mapsto\mathbb{S}^{d+1}_{\succeq0}$,  $R:[0,T]\mapsto\mathbb{S}^{m}_{\succ0}$, and $Q_T\in\mathbb{S}^{d+1}_{\succeq0}$.
%
%
We then have the following result:
\begin{theorem}\label{th:uni:CT:LQ}
  Let us consider the unimolecular reaction network $(\X{},\mathcal{R})$ with propensity functions \eqref{eq:propensities:unimolecular:inputs}. Then, the optimal control input minimizing the cost  \eqref{eq:LQcost} is given by
  \begin{equation}\label{eq:uopt:CT}
  u^*(t)=-\dfrac{1}{2}R(t)^{-1}\Omega_u(t)\bar x(t)
  %
%
\end{equation}
where 
\begin{equation}
  \begin{array}{rcl}
   \Omega_u(t)&:=&2\bar B(t)^{\T}P(t)+\Gamma(t)e_{d+1}^{\T},\\
   \Gamma(t)&:=&\sum_{k=1}^K\bar{B}_k(t)^{\T}P(t)\bar\zeta_k,
  \end{array}
\end{equation}
 and $P:[0,T]\mapsto\mathbb{S}^{d+1}$, $P(T)=Q_T$, solves the Riccati differential equation
\begin{equation}\label{eq:Riccati:uniCT}
\begin{array}{l}
  \dot{P}(t)+\bar A(t)^{\T}P(t)+P(t)\bar A(t)+Q(t) -\dfrac{1}{4}\Omega_u(t)^{\T}R(t)^{-1}\Omega_u(t)\\
  \qquad+\dfrac{1}{2}\sum_{k=1}^K(\bar A_k(t)^{\T}P(t)\bar\zeta_ke_{d+1}^{\T}+e_{d+1}\bar\zeta_k^{\T}P(t)\bar A_k(t))=0.
\end{array}
\end{equation}
where $\bar A_k(t)=\bar\zeta_k \bar w_{k}(t)^{\T}$, $\bar B_k(t)=\bar\zeta_k \bar b_{k}(t)^{\T}$, 
\begin{equation}
  \bar A(t)=\begin{bmatrix}
    A(t) & d(t)\\
    0 & 0
  \end{bmatrix},\ \bar B(t)=\begin{bmatrix}
    B(t)\\
    0
  \end{bmatrix},\textnormal{ and }\bar\zeta_k=\begin{bmatrix}
    \zeta_k\\
    0
  \end{bmatrix}.
\end{equation}
Moreover, we have that $J_T^*=J_T(u^*)=V(0,X(0))$.
\end{theorem}
\begin{proof}
Since the running cost is strictly convex in $u$ and the propensity functions are affine in $u$, then a minimum for $J_T(u)$ exists. We will show that a quadratic value function of the form $V(t,x)=\bar x^{\T}P(t)\bar x$ verifies the \HJB equation \eqref{eq:HBJCT} which is given in this case by
  \begin{equation}\label{eq:HBJCT_uni}
  \bar x^{\T}\dot{P}(t)\bar x +\bar x^{\T}Q(t)\bar x+\min_u\left\{\sum_{i=1}^K\lambda_i(t,x,u)\Delta_iV(t,x)+u^{\T}R(t)u\right\}=0
\end{equation}
and $P(T)=Q_T$ where
\begin{equation}
    \Delta_iV(t,x)=\bar\zeta_i^{\T}P(t)\bar x+\bar x^{\T}P(t)\bar\zeta_i+\bar\zeta_i^{\T}P(t)\bar\zeta_i
\end{equation}
and where $ \bar\zeta_k$ is defined in the result.
Computing the derivative with respect to $u$ and equating it to zero yields
\begin{equation}
  \sum_{i=1}^K\Delta_iV(t,x)\bar{b}_i(t)+2R(t)u=0,
\end{equation}
which results in the optimal control law
\begin{equation}
  u^*(t,x):=\dfrac{1}{2}R(t)^{-1}\left(\sum_{i=1}^K\Delta_iV(t,x)\bar{b}_i(t)\right).
\end{equation}
Substituting the explicit expressions for the $ \Delta_iV(t,x)$'s yields \eqref{eq:uopt:CT} and substituting, in turn, \eqref{eq:uopt:CT} in \eqref{eq:HBJCT_uni} yields
\begin{equation}
\begin{array}{l}
   \bar x^{\T}\dot{P}(t)\bar x-\dfrac{1}{4}(2\bar B(t)^{\T}P(t)\bar x+\Gamma(t))^{\T}R(t)^{-1}(2\bar B(t)^{\T}P(t)\bar x+\Gamma(t))\\
  +\sum_{k=1}^K\bar w_i(t)\bar x\left(2\bar\zeta_k^{\T}P(t)\bar x+\bar\zeta_k^{\T}P(t)\bar\zeta_k\right)+\bar{x}^{\T}Q(t)\bar{x}=0.
\end{array}
\end{equation}
Using the fact that
\begin{equation}
  \begin{array}{rcl}
  \bar w_k(t)^{\T}\bar x\left(2\bar\zeta_k^{\T}P\bar x+\bar\zeta_k^{\T}P\bar\zeta_k\right) &=&\bar x^{\T}\left(\He\left[P(t)\bar\zeta_k\bar w_k(t)^{\T}+\dfrac{1}{2}\bar w_k(t)\bar\zeta_k^{\T}P(t)\bar\zeta_ke_{d+1}^{\T}\right]\right)\bar x\\
    &=&\bar x^{\T}\left(\He\left[P(t)\bar A_k(t)+\dfrac{1}{2}\bar A_k(t)^{\T}P(t)\bar{\zeta}_ke_{d+1}^{\T}\right]\right)\bar x,
  \end{array}
\end{equation}
we obtain
\begin{equation}
  \begin{array}{rcl}
  \sum_{k=1}^K \lambda_k(t,x)\left(2\bar\zeta_k^{\T}P(t) \bar x+\bar\zeta_k^{\T}P(t)\bar\zeta_k\right) &=&\bar x^{\T}\left(\He\left[P(t)\bar A(t)+\dfrac{1}{2}\sum_{k=1}^K\bar A_k(t)^{\T}P(t)\bar{\zeta}_ke_{d+1}^{\T}\right]\right)\bar x.
  \end{array}
\end{equation}
Finally, noting that $\Gamma(t)=\Gamma(t) e_{d+1}^{\T}\bar x$ yields
\begin{equation}
\begin{array}{l}
   \bar x^{\T}\dot{P}(t)\bar x-\bar x^{\T}\dfrac{1}{4}(2\bar B(t)^{\T}P(t)+\Gamma(t) e_{d+1}^{\T})^{\T}R(t)^{-1}(2\bar B(t)^{\T}P(t)+\Gamma(t) e_{d+1}^{\T})\bar x\\
  +\bar x^{\T}\left(\He\left[P(t)\bar A(t)+\dfrac{1}{2}\sum_{k=1}^K\bar A_k(t)^{\T}P(t)\bar{\zeta}_ke_{d+1}^{\T}\right]\right)\bar x+\bar{x}^{\T}Q(t)\bar{x}=0.
\end{array}
\end{equation}
As this expression holds for all $\bar{x}\in\mathbb{Z}_{\ge0}^d$, then this is equivalent to the Riccati differential equation \eqref{eq:Riccati:uniCT}.
\end{proof}

As a matter of comparison, we can state here the more classical result on the optimal LQ control of the moment equation \eqref{eq:moment:uni:CT}:
\begin{theorem}\label{th:ELQ}
  Let us consider the moment system \eqref{eq:moment:uni:CT}. Then, the optimal (deterministic) control input that minimizes the cost
  \begin{equation}
    J_T^E(u):=\int_0^T\left(\E[\bar X(s)]^{\T}Q(s)\E[\bar X(s)]+u(s)^{\T}R(s)u(s)\right)\ds+\E[\bar X(T)]^{\T}Q_T\E[\bar X(T)]
  \end{equation}
  is given by
  \begin{equation}
    u^*(t):=-R(t)^{-1}\bar B(t)^{\T}P(t)\E[\bar X(t)]
  \end{equation}
  where $P:[0,T]\mapsto\mathbb{S}^{d+1}$, $P(T)=Q_T$, solves the Riccati differential equation
  \begin{equation}
    \dot{P}(t)+\bar A(t)^{\T}P(t)+P(t)\bar A(t)+Q(t)-P(t)\bar B(t)R(t)^{-1}\bar B(t)P(t)=0,\ t\in[0,T].
  \end{equation}
  Moreover, we have that $J_T^{E*}:=J_T^E(u^*)=\E[\bar{X}_0]^{\T}P(0)\E[\bar{X}_0]$.
\end{theorem}

It is important to stress that the obtained Riccati differential equation obtained in Theorem \ref{th:uni:CT:LQ} is dramatically different from the classical one provided in Theorem \ref{th:ELQ} and those obtained in the LQ control of linear stochastic differential equations. Firstly, we can observe the presence of the additional $\Gamma(t)$ terms in both the optimal control law and in the Riccati differential equation. Secondly, there is an extra sum, linear in $P(t)$, in the Riccati differential equation. This is a consequence of the fact that jumps act on the state in the additive way $x \mapsto x+\zeta_k$.

It seems also important to clarify under which conditions a solution to the Riccati differential equation exists on the interval $[0,T]$ as such equations may have a finite escape time. From the theory of Riccati differential equations \cite{AbouKandil:03}, we know that a solution will exist in an interval $(T-\eps,T]$ for some sufficiently small $\eps>0$. However, using the fact that the matrix $P$ is symmetric and positive semidefinite, the standard arguments for the global existence of a solution $P$ on $[0,T]$ can be straightforwardly adapted to the current Riccati equation. However, it seems that the characterization of the solution in terms of a Hamiltonian matrix does not exist in this setting due to the very peculiar structure of the Riccati differential equation.}

\section{Optimal Sampled-Data Control of Stochastic Reaction Networks}\label{sec:SD:LQ}

\blue{There are multiple ways for considering sampled-data systems and sampled-data control problems. Chronologically, the first approach to have been introduced was based on the direct discretization of the dynamics of the process \cite{Jury:58}. While this procedure is quite straightforward in the linear deterministic case, it is way more involved in the stochastic or the nonlinear settings because of their more complex nature. A subsequently proposed approach was based on the so-called "\textit{lifting procedure}" which turned the sampled-data system into an infinite-dimensional system \cite{Yamamoto:90,Yamamoto:94}. This approach had the benefit of preserving the information on the \textit{inter-sample  behavior}, lost by the discretization procedure, at the expense of having to work in the infinite-dimensional setting. Variations of the approach managed to keep the dimension of the state-space finite but the input and output spaces still had to be considered as infinite-dimensional \cite{Bamieh:91,Chen:95} in order to be able to properly address the performance analysis and optimization problems. Another approach that makes the system infinite-dimensional is the so-called \textit{input-delay approach} that turns the sampled-data system into a system with a sawtooth time-varying delay \cite{Fridman:04} and on which available tools for the analysis and the control of time-delay systems can be applied. An impulsive systems formulation was considered in \cite{Sivashankar:94} and had the advantage of keeping the state-space as well as the input- and output spaces finite-dimensional, at the reasonably affordable expense of making the dynamics of the system hybrid. This framework is now included  in the very general hybrid systems framework \cite{Goebel:12} for which many analysis tools have been developed since then.

The main limitation for using a discrete-time approach in the current setting lies in the difficulty in computing or even considering the map  $F:\mathbb{R}_{\ge0}\times\mathbb{Z}_{\ge0}^d\times\mathbb{R}^m_{\ge0}\mapsto\mathbb{Z}_{\ge0}^d$ such that
\begin{equation}
    X(t_{k+1})=F(k,X(t_k),u(t_k)),\ t_k=kT_s
\end{equation}
 where $T_s>0$ is the sampling period and $u(t_k)$ is the control input value at $t_k$. This map clearly relates to the transition kernel $\P(X(t_{k+1})=x^+|X(t_k)=x)$ which is, in general, a complex function of the propensity functions related to the exponential of the CME operator which is not exactly computable in the current setting. Even if it were computable, the resulting dependence on $u(t_k)$ would be quite complex and, hence, difficult to consider for optimal control purposes. Finally, the inter-sample behavior would be lost when discretizing the process while this information is necessary to consider in the cost to minimize. The lifting and the input-delay approaches are much less convenient to use than the hybrid formulation which keeps the state-space finite-dimensional and allows for the use of simpler concepts than those required for dealing with infinite-dimensional systems. Adding the fact that a Markov jump process can already be seen as a stochastic hybrid system closes the deal in favor of the hybrid systems approach.}

\subsection{Sampled-data stochastic reaction networks as stochastic hybrid systems}

\blue{A sampled-data stochastic reaction network is stochastic reaction network having a piecewise constant control input which is only updated at every time $t_k=kT_s$, $k\in\mathbb{Z}_{\ge0}$, where $T_s>0$ is the sampling period. Such a system can be described as a stochastic hybrid system \cite{Teel:14} with
\begin{itemize}
  \item a \textbf{flow map}
\begin{equation}\label{eq:SD:1}
\left.\begin{array}{rcl}
  \dot{X}(t)&=&0\\
  \dot{v}(t)&=&0
\end{array}\right\} \textnormal{ if }t\in[0,T]-\{t_k\}_{k\ge0}
\end{equation}
\item \textbf{$K$ spontaneous jump maps}
\begin{equation}\label{eq:SD:2}
\left.\begin{array}{rcl}
  X(t^+)&=&X(t)+\zeta_i\\
  v(t^+)&=&v(t)
  \end{array}\right\}\sim \lambda_i(t,X(t),v(t)),\ i=1,\ldots,K,\ \textnormal{and}
\end{equation}
\item a \textbf{deterministic jump map}
\begin{equation}\label{eq:SD:3}
\left.\begin{array}{rcl}
  X(t^+)&=&X(t)\\
  v(t^+)&=&u(t)
  \end{array}\right\} \textnormal{ if }t\in \{t_k\}_{k\ge0}\cap[0,T].
\end{equation}
\end{itemize}
%
%
It seems important to explain this model in more details. Since the molecular counts are piecewise constant, it is immediate to see that the derivative of the molecular counts between jumps should be zero. The second state, $v\in\mathbb{R}_{\ge0}^d$, has been added in order to model the zero-order hold (ZOH) that keeps the value of the control input constant over each sampling period. In this regard, the flow map can be considered as trivial in the current context. The $K$ spontaneous jumps, which are exponentially distributed with rate $\lambda_i(t,X,v)$, $i=1,\ldots,K$, correspond to the random firings of the $K$ possible reactions and only change the value of the molecular counts. Finally, the deterministic jump map updates the value of the control input while leaving the molecular counts unchanged.}

\subsection{General stochastic reaction networks}

\blue{In the case of general networks, the propensity functions can be written as $\lambda_i(t,x,v)$ as defined in \eqref{eq:propensities_bimolecular} where we replace the control input $u$ the state $v$ of the ZOH. We consider here the following cost
\begin{equation}\label{eq:hybridcostSRN}
  J_T(u):=\E\left[\int_{0}^{T}q_{c}(s,x(s),v(s))\ds+\sum_{k=0}^{N-1}q_{d}(k,x(t_k),v(t_k),u(t_k))+q_{c,T}(x(T),v(T))\right]
\end{equation}
where $q_c:[0,T]\times \mathbb{Z}_{\ge0}^d\times \mathbb{R}_{\ge0}^m\mapsto\mathbb{R}_{\ge0}$ is the continuous-time running cost, $q_d:\{0,\ldots,N-1\}\times \mathbb{Z}_{\ge0}^d\times \mathbb{R}_{\ge0}^m\times \mathbb{R}_{\ge0}^m\mapsto\mathbb{R}_{\ge0}$ is the discrete-time running cost, and $q_{c,T}: \mathbb{Z}_{\ge0}^d\times \mathbb{R}_{\ge0}^m\mapsto\mathbb{R}_{\ge0}$ is the terminal cost.  We also assume, without loss of generality, that $t_N=NT_s=T$. The control law is then given by
\begin{equation}
  u^*=\argmin_{u\in\mathcal{V}_{0,T}} J_T(u)
\end{equation}
where we have assumed that such a minimum exists and
\begin{equation}
\mathcal{V}_{s,t}:=\left\{u\in \ell_2([s,t]\cap\{t_0,\ldots,t_N\}\times\mathbb{Z}_{\ge0}^d,\mathbb{R}_{\ge0}^m):\ u\ \textnormal{mesurable}\right\},\ t\ge s.
\end{equation}
Considering the value function
\begin{equation}
\begin{array}{rcl}
  V(t,x,v) &=& \min_{u\in\mathcal{V}_{t,T}}\E\left[\left.\int_t^{T}q_c(s,x(s),v(s))\ds+\smashoperator{\sum_{k:t<t_k\le (N-1)T}}q_{d}(k,x(t_k),v(t_k),u(t_k))\right|\begin{array}{rcl}
    X(t)&=&x\\
    v(t)&=&x
  \end{array}\right]
\end{array}
\end{equation}
which is defined for all $(x,v)\in\mathbb{Z}_{\ge0}^d\times\mathbb{R}_{\ge0}^m$, we obtain the following result:

\begin{theorem}\label{th:generalSD}
Let us consider the sampled-data reaction network represented by  \eqref{eq:SD:1}, \eqref{eq:SD:2}, \eqref{eq:SD:3}, and \eqref{eq:propensities_bimolecular}. The optimal control input is given by
\begin{equation}
  u^*(t_k)=\argmin_u \{V(t_k^+,x,u)+q_{d}(k,x,u)\},\ k=0,\ldots,N-1
\end{equation}
where the value function $V$ verifies the following hybrid HJB equation
\begin{equation}
\begin{array}{l}
    V_t(t,x,v)+\sum_{k=1}^K\lambda_k(t,x,v)(V(t,x+\zeta_k,v)-V(t,x,v))\\
  \qquad\qquad+q_{c}(t,x,v)=0,\ t\in(t_k,t_{k+1}],\ k=0,\ldots,N-1
\end{array}
\end{equation}
and
\begin{equation}
  V(t_k^+,x,u^*(t_k))-V(t_k,x,v)+q_{d}(k,x,v,u^*(t_k))=0,\ k=0,\ldots,N-1
\end{equation}
with the terminal condition $V(T,x,v)=q_{c,T}(x,v)$.
\end{theorem}
\begin{proof}
  The proof simply consists of applying Dynamic Programming to the system  \eqref{eq:SD:1}, \eqref{eq:SD:2}, \eqref{eq:SD:3}, and \eqref{eq:propensities_bimolecular}.
\end{proof}}

\subsection{Unimolecular stochastic reaction networks}

\blue{In this case, the propensity functions vector is given by
\begin{equation}\label{eq:propensities_unimolecular_hybrid}
  \lambda_i(t,x,v)=\bar w_i(t)^{\T}\bar x+\bar b_i(t)^{\T}v=:\bar{\lambda}_i(t)^{\T}z
\end{equation}
where $\bar w_i(t) = \col(w_{i,x}(t),w_{i,0}(t))$, $w_{i,x}\in\mathbb{R}_{\ge0}^{d}$, $\bar w_{i,0}(t)\in\mathbb{R}_{\ge0}$, $\bar b_i(\cdot)\in\mathbb{R}_{\ge0}^{m}$,  $\bar{\lambda}_i(t):=\col(\bar w_i(t),\bar b_i(t))$, $\bar x:=\col(x,1)$, and $z:=\col(x,1,v)$. As in the continuous-time case, we can write the moment equations in closed-form with the difference that it takes now the form of a linear impulsive system
\begin{equation}
\begin{array}{rcl}
  \dfrac{\d}{\dt}\begin{bmatrix}
    \E[X(t)]\\
    \E[v(t)]
  \end{bmatrix}&=&\begin{bmatrix}
    A(t) & B(t)\\
    0 & 0
  \end{bmatrix}\begin{bmatrix}
    \E[X(t)]\\
    \E[v(t)]
  \end{bmatrix}+\begin{bmatrix}
    d(t)\\
    0
  \end{bmatrix},\ t\ne t_k\\
  \begin{bmatrix}
    \E[X(t_k^+)]\\
    \E[v(t_k^+)]
  \end{bmatrix}&=&\begin{bmatrix}
    I & 0\\
    0 & 0
  \end{bmatrix}\begin{bmatrix}
    \E[X(t_k)]\\
    \E[v(t_k)]
  \end{bmatrix}+\begin{bmatrix}
    0\\
    I
  \end{bmatrix}\E[u(t_k)],\ k=0,\ldots,N-1
\end{array}
\end{equation}
where
\begin{equation}
  A(t)=\sum_{k=1}^K\zeta_k \bar w_{k,x}(t)^{\T},\ B(t)=\sum_{k=1}^K\zeta_k \bar b_k(t)^{\T},\ \textnormal{and }  d(t)=\sum_{k=1}^K\zeta_k \bar w_{k,0}(t)^{\T}.
\end{equation}

We also introduce the following hybrid LQ cost
\begin{equation}\label{eq:LQcostHybrid}
\begin{array}{rcl}
    J_T(u)&:=&\E\left[\int_{0}^{T}z(s)^{\T}Q_{c}(s)z(s)\ds+\sum_{k=0}^{N-1}\left(z(t_k)^{\T}Q_{d}(k)z(t_k)+u(t_k)^{\T}R_{d}(k)u(t_k)\right)\right.\\
    &&\left.\vphantom{\int_{0}^{T}z(s)^{\T}Q_{c}(s)z(s)\ds}\qquad+z(T)^{\T}Q_{c,T}z(T)\right].
\end{array}
\end{equation}
where $Q_{c}:[0,T]\mapsto\mathbb{S}^{d+1+m}_{\succeq0}$, $Q_{c,T}\in\mathbb{S}^{d+1+m}_{\succeq0}$, $Q_{d}:\{0,\ldots,N-1\}\mapsto\mathbb{S}^{d+1+m}_{\succeq0}$ and $R_d:\{0,\ldots,N-1\}\mapsto\mathbb{S}^{m}_{\succ0}$. This leads to the following result:
\begin{theorem}\label{th:SDuni}
Let us consider the sampled-data unimolecular stochastic reaction network \eqref{eq:SD:1},  \eqref{eq:SD:2}, and \eqref{eq:SD:3} with propensities \eqref{eq:propensities_unimolecular_hybrid}. Then, the optimal control input that minimizes the cost \eqref{eq:LQcostHybrid} is given by
  \begin{equation}\label{eq:uopt:DT}
  u^*(t_k)=K(k) z(t_k)=-(R_d(k)+P_{3}(t_k^+))^{-1} P_2(t_k^+)^{\T}\bar x(t_k),\ k=0,\ldots,N-1
\end{equation}
where $P:[0,T]\mapsto\mathbb{S}^{d+1+m}$, $P(T)=Q_T$, solves the Lyapunov differential equation
\begin{equation}\label{eq:Riccati:uniDT1}
\begin{array}{l}
  \dot{P}(t)+Q_{c}(t)+\bar{A}(t)^{\T}P(t)+P(t)\bar{A}(t)+\dfrac{1}{2}\sum_{i=1}^K\Sym(\bar{\lambda}_i(t)\bar\zeta_i^{\T}P(t)\bar\zeta_ie_{d+m+1}^{\T})=0, t\in(t_k,t_{k+1}]
\end{array}
\end{equation}
and the Riccati difference equation
\begin{equation}\label{eq:Riccati:uniDT2}
  \begin{array}{rcl}
    P_{1}(t_k^+)-P_{2}(t_k^+)(P_3(t_k^+)+R_d(k))^{-1}P_2(t_k^+)^{\T}+Q_{d,1}(k)-P_{1}(t_k)&=&0\\
    Q_{d,2}(k)-P_{2}(t_k)&=&0\\
    Q_{d,3}(k)-P_{3}(t_k)&=&0
  \end{array}
\end{equation}
for $k=0,\ldots,N-1$ and where $\bar{\zeta}_i=\col(\zeta_i,0,0)$ together with
\begin{equation}\label{eq:Partitiojnr}
  \bar A(t)=\begin{bmatrix}
    A(t) & d(t) & \vline & B(t)\\
    0 & 0 &\vline & 0\\
    \hline
    0 & 0 & \vline & 0
  \end{bmatrix},\ P(t)=\begin{bmatrix}
    P_1(t) & \vline & P_2(t)\\
    \hline
    P_2(t)^{\T} & \vline &P_3(t)
  \end{bmatrix},\ Q_d(k) = \begin{bmatrix}
     Q_{d,1}(k) &  \vline & Q_{d,2}(k)\\
     \hline
     Q_{d,2}(k)^{\T} &  \vline &Q_{d,3}(k)
  \end{bmatrix}
\end{equation}
with $P_1(t)\in\mathbb{S}^{d+1}$, $P_2(t)\in\mathbb{R}^{(d+1)\times m}$, $P_3(t)\in\mathbb{S}^{m}$, $Q_{d,1}(k)\in\mathbb{S}_{\succeq0}^{d+1}$, $Q_{d,2}(k)\in\mathbb{R}^{(d+1)\times m}$, and $Q_{d,3}(k)\in\mathbb{S}_{\succeq0}^{m}$.
\end{theorem}
\begin{proof}
 Let $\bar\zeta_i=(\zeta_i,0,0)\in\mathbb{R}^{d+m+1}$. We will show that the value function $V(t,x,v):=z^{\T}P(t)z$, $P(t)\in\mathbb{S}^{d+1+m}_{\succeq0}$, verifies the hybrid Dynamic Programming equations of Theorem \ref{th:generalSD}. For the flow part, we get that
\begin{equation}
  z^{\T}\dot{P}(t)z+z^{\T}Q_{c,k}(t)z+\sum_{i=1}^K\lambda_i(t,x,v)\Delta_iV(t,x,v)=0
\end{equation}
where
\begin{equation}
  \Delta_iV(t,x,v):=\bar\zeta_i^{\T}P(t)z+\bar z^{\T}P(t)\bar\zeta_i+\bar\zeta_i^{\T}P(t)\bar\zeta_i.
\end{equation}
Therefore
\begin{equation}\label{eq:kdslkd;sk;dkkd;a;sdk}
  z^{\T}\dot{P}(t)z+z^{\T}Q_{c}(t)z+\sum_{i=1}^K  \bar{\lambda}_i(t)^{\T}z\left(\bar\zeta_i^{\T}P(t)z+z^{\T}P(t)\bar\zeta_i+\bar\zeta_i^{\T}P(t)\bar\zeta_i\right)=0.
\end{equation}
The above expression can be decomposed as the sum of
\begin{equation}
  \sum_{i=1}^K\bar{\lambda}_i(t)^{\T}z\left(\bar\zeta_i^{\T}P(t)z+z^{\T}P(t)\bar\zeta_i\right)=\bar A(t)^{\T}P(t)+P(t)\bar A(t)
\end{equation}
with
\begin{equation}
  \sum_{i=1}^K\bar{\lambda}_i(t)^{\T}z\bar\zeta_i^{\T}P(t)\bar\zeta_i=\dfrac{1}{2}\sum_{i=1}^Kz^{\T}\Sym(\bar{\lambda}_i(t)\bar\zeta_i^{\T}P(t)\bar\zeta_ie_{d+1}^{\T})z.
\end{equation}
Substituting those two expressions in \eqref{eq:kdslkd;sk;dkkd;a;sdk} yields the Lyapunov differential equation \eqref{eq:Riccati:uniDT1}. This proves the result for the flow part.\\


We now look at the jump part of the dynamics. The Dynamic Programming equation is given in this case by
\begin{equation}\label{eq:dksl;kdl;skd;kdl;sdl;ksa}
 \min_{u(t_k)} \begin{bmatrix}
    x(t_k)\\
    1\\
    u(t_k)
  \end{bmatrix}^{\T}P(t_k^+)\begin{bmatrix}
    x(t_k)\\
    1\\
    u(t_k)
  \end{bmatrix}-z(t_k)^{\T}P(t_k)z(t_k)+z(t_k)^{\T}Q_d(k)z(t_k)+u(t_k)^{\T}R_d(k)u(t_k)=0.
\end{equation}
Let us then partition $P$ as in \eqref{eq:Partitiojnr} and computing the derivative with respect to $u(t_k)$ yields
\begin{equation}
\begin{array}{rcl}
  \dfrac{\partial}{\partial u(t_k)} \left(\begin{bmatrix}
    x(t_k)\\
    1\\
    u(t_k)
  \end{bmatrix}^{\T}P(t_k^+)\begin{bmatrix}
    x(t_k)\\
    1\\
    u(t_k)
  \end{bmatrix}+u(t_k)^{\T}R_d(k)u(t_k)\right)&=&2\begin{bmatrix}
    P_{2}(t_k^+)^{\T} & \vline & P_{3}(t_k^+)
  \end{bmatrix}\begin{bmatrix}
    x(t_k)\\
    1\\
    \hline
    u(t_k)
  \end{bmatrix}\\

  &&\qquad+2u(t_k)^{\T}R_d(k).
\end{array}
\end{equation}
This means that the optimal control input verifies the expression
\begin{equation}
  (R_d(k)+P_{3}(t_k^+))u(t_k)+P_{2}(t_k^+)^{\T}\bar x(t_k)=0
\end{equation}
which yields, in turn, the expression $u(t_k)=K(k)\bar x(t_k)$ with $K(k)=-(R_d(k)+P_{3}(t_k^+))^{-1}P_{2}(t_k^+)^{\T}$. If we now define $H(k)$ as
\begin{equation}
  H(k):=\begin{bmatrix}
    I_{d+1} & \vline & 0\\
    \hline
     K(k)&\vline& 0
  \end{bmatrix}
\end{equation}
then, we get that
\begin{equation}
  \begin{bmatrix}
    x(t_k)\\
    1\\
    \hline
    u(t_k)
  \end{bmatrix}=H(k)z(k),
\end{equation}
and, substituting this expression in \eqref{eq:dksl;kdl;skd;kdl;sdl;ksa}, we obtain the expression
\begin{equation}\label{eq:jdlsjdsldjal}
 H(k)^{\T}P(t_k^+)H(k)+Q_d(k)+\begin{bmatrix}
   K(k) & 0
 \end{bmatrix}^{\T}R_d(k)\begin{bmatrix}
   K(k) & 0
 \end{bmatrix}-P(t_k)=0.
\end{equation}

Expanding the above equality yields
\begin{equation}
  \begin{bmatrix}
    P_{1}(t_k^+)+\He[P_{2}(t_k^+)K(k)]+ K(k)^{\T}(R_d(k)+P_{3}(t_k^+)K(k) & 0\\
    0 & 0
  \end{bmatrix}+Q_d(k)-P(t_k)=0,
\end{equation}
which reduces to $Q_{d,2}(k)- P_{2}(t_k)=0$, $Q_{d,3}(k)-  P_{3}(t_k)=0$ , and
\begin{equation}
     P_{1}(t_k^+)+\He[P_{2}(t_k^+)K(k)]+ K(k)^{\T}(R_d(k)+P_{3}(t_k^+)K(k)+Q_{d,1}(k)-  P_{1}(t_k)=0.
\end{equation}
Finally, since $K(k)=-(R_d(k)+P_{3}(t_k^+))^{-1}P_{2}(t_k^+)^{\T}$, the above expression simplifies to
\begin{equation}
    P_{1}(t_k^+)-P_{2}(t_k^+)(P_3(t_k^+)+R_d(k))^{-1}P_2(t_k^+)^{\T}+Q_{d,1}(k)-  P_{1}(t_k)=0
\end{equation}
together with $Q_{d,2}(k)-  P_{2}(t_k)=0$, $Q_{d,3}(k)-  P_{3}(t_k)=0$, which proves \eqref{eq:Riccati:uniDT2} and the desired result.
\end{proof}}

\section{Continuous-Time $H_\infty$ Control of Stochastic Reaction Networks}\label{sec:CT:Hinf}

The objective of this section is to develop an  $H_\infty$ (or $L_2$) control theory for stochastic reaction networks. The key idea is to derive a control law that can guarantee a certain attenuation level from exogenous inputs to control outputs in an $L_2$-sense. Such a result exists in a wide variety of contexts (see e.g. \cite{vanderSchaft:00}) but does not seem to have been developed for jump Markov processes describing stochastic reaction networks. A general result is first derived and is specialized to the unimolecular case, for which a Bounded Real Lemma is also obtained.

\subsection{Preliminaries}

\blue{Similarly as in the previous section, we extend the reaction network \net to not only involve $m$ control inputs but also $p$ exogenous inputs, all of them acting at the propensity level. This gives rise to the following structural assumption:
\begin{hyp}\label{hyp:CTHinf1}
  The propensity functions of the network \net are of the form $ \lambda(t,x,u,w)$ where $u\in\mathbb{R}^m_{\ge0}$ is the control input vector and $w\in\mathbb{R}_{\ge0}^p$ is the vector of exogenous disturbances. We also define the so-called controlled output vector  $z(t)=h(t,x(t),u(t),w(t))$ where $h:\mathbb{R}_{\ge0}\times\mathbb{Z}_{\ge0}^d\times\mathbb{R}^m_{\ge0}\times\mathbb{R}_{\ge0}^p\mapsto\mathbb{R}^q$ is assumed to be a continuous function.
\end{hyp}
%
%
%

The $H_\infty$ control is tightly connected to the concept of $L_2$-gain of operators and the $L_2$-norm of signals. Before being able to proceed with the main results, it is important to state the following definitions:
\begin{define}
  The (finite-time) $L_2$-norm of the signal $w:[0,T]\mapsto\mathbb{R}^p$ is defined as
\begin{equation}
  ||w||_{L_2([0,T])}:=\left(\int_0^T\E[||w(s)||_2^2]\ds\right)^{1/2}.
\end{equation}
\end{define}

We can now define the concept of $L_2$-gain of operators:
\begin{define}
The finite-horizon  $L_2$-gain of the transfer/operator $w\mapsto z$ associated with the stochastic reaction network \net under Assumption \ref{hyp:CTHinf1} is defined as
\begin{equation}
  ||w\mapsto z||_{L_2([0,T])-L_2([0,T])}:=\sup_{w\in \mathcal{W}_{0,T}, ||w||_{L_2([0,T])}\ne0,X(0)=0}\dfrac{||z||_{L_2([0,T])}}{||w||_{L_2([0,T])}}
\end{equation}
where
\begin{equation}
  \mathcal{W}_{s,t}:=\left\{w\in L_2([s,t]\times\mathbb{Z}_{\ge0}^d\times\mathbb{R}^m_{\ge0},\mathbb{R}_{\ge0}^p):\ w\ \textnormal{mesurable}\right\},\ t\ge s.
\end{equation}
\end{define}

The $H_\infty$ control problem consists of solving the following optimization problem
\begin{equation}\label{eq:Hinfcontrolproblem}
  \min_{u\in\mathcal{U}_{0,T}}||w\mapsto z||_{L_2([0,T])-L_2([0,T])}=\min_{u\in\mathcal{U}_{0,T}}\sup_{w\in \mathcal{W}_{0,T}, ||w||_{L_2([0,T])}\ne0,X(0)=0}\dfrac{||z||_{L_2([0,T])}}{||w||_{L_2([0,T])}},
\end{equation}
where $\mathcal{U}_{s,t}$ is defined in \eqref{eq:Ust}. It is important to stress that in the above optimal control problem, the control input is not allowed to depend on the (typically unknown) exogenous input $w$, while the reverse is possible.

\subsection{General stochastic reaction networks}

A way to solve the $H_\infty$ control problem \eqref{eq:Hinfcontrolproblem} is through the introduction of the cost
\begin{equation}\label{eq:generalHinfcost}
\begin{array}{rcl}
  J_T(u,w)&:=&\E\left[\displaystyle\int_0^T\left(r(s,u(s))+||z(s)||_2^2-\gamma^2||w(s)||_2^2\right)\ds+q_T(X(T))\right]
\end{array}
\end{equation}
where $r:[0,T]\times \mathbb{R}_{\ge0}^m\mapsto\mathbb{R}_{\ge0}$ is the running cost for the control input and $q_T:\mathbb{Z}_{\ge0}^d\mapsto\mathbb{R}$ is the terminal cost. The control input is chosen in such a way that it minimizes the worst case
\begin{equation}\label{eq:Hinfminimax}
  J_T^*=\min_{u\in\mathcal{U}_{0,T}}\max_{w\in\mathcal{W}_{0,T}}  J_T(u,w)
\end{equation}

We also make the following assumption
\begin{hyp}\label{hyp:CTHinf2}
The equality
  \begin{equation}
  \dfrac{\partial^2}{\partial u\partial w}\left(\lambda(t,x,u,w)^{\T}\Delta V(t,x)+r(t,u)+z^{\T}z-\gamma^2w^{\T}w\right)=0
\end{equation}
  holds and the function $\lambda(t,x,u,w)^{\T}\Delta V(t,x)+r(t,u)+z^{\T}z-\gamma^2w^{\T}w$ is strongly convex in $u$ and strongly concave in $w$.
\end{hyp}

Considering the value function
\begin{equation}
  V(t,x) = \min_{u\in\mathcal{U}_{t,T}}\max_{w\in\mathcal{W}_{t,T}}\E\left[\left.\int_t^T\left(r(s,u(s))+||z(s)||_2^2-\gamma^2||w(s)||_2^2\right)\ds\right|X(t)=x\right],\ V(T,x)=q_T(x)
\end{equation}
defined for all $x\in\mathbb{Z}_{\ge0}^d$ leads to the following result:
\begin{theorem}
Let us consider a stochastic reaction network satisfying Assumption \ref{hyp:CTHinf1} and Assumption \ref{hyp:CTHinf2}. Then an optimal pair $(u^*,w^*)$ that solves the optimization problem \eqref{eq:Hinfminimax} is given by
\begin{equation}\label{eq:Hinfoptimalpoint}
\begin{array}{rcl}
  u^*(t)&=&\argmin_u\left\{\sum_{i=1}^K\lambda_i(t,x(t),u,w)\Delta_iV(t,x(t))+r(t,u)+z(t)^{\T}z(t)\right\},\\
  w^*(t)&=&\argmax_w\left\{\sum_{i=1}^K\lambda_i(t,x(t),u,w)\Delta_iV(t,x(t))+z(t)^{\T}z(t)-\gamma^2w^{\T}w\right\},
\end{array}
\end{equation}
where the function $V(t,x)$ satisfies the Hamilton-Jacobi-Isaacs (HJI) equation given by
\begin{equation}\label{eq:HBICT}
  V_t(t,x)+\min_{u}\max_{w} \left\{\sum_{i=1}^K\lambda_i(t,x,u,w)\Delta_iV(t,x)+r(t,u)+z(t)^{\T}z(t)-\gamma^2w(t)^{\T}w(t)\right\}=0
\end{equation}
and $V(T,x)=q_T(x)$.

Moreover, in such a case, the optimal cost $J_T^*:=J_T(u^*,w^*)$ coincides with $V(0,X(0))$ and we have that
\begin{equation}
 V(0,X(0))= ||z||_{L_2([0,T])}^2-\gamma^2||w^*||_{L_2([0,T])}^2+\int_0^Tr(s,u^*(s))\ds+\E[q_T(X(T))].
\end{equation}
which reduces in the case of no terminal cost and zero initial condition to
\begin{equation}
 ||z||_{L_2([0,T])}^2<\gamma^2||w||_{L_2([0,T])}^2
\end{equation}
for all $w\in \mathcal{W}_{0,T}$.
\end{theorem}
\begin{proof}
  Due to Assumption \ref{hyp:CTHinf2} there exists an unique optimal point $(u^*,w^*)$ given by \eqref{eq:Hinfoptimalpoint} to the optimal control problem \eqref{eq:Hinfminimax} and the optimal control input $u^*$ is independent of $w$. The rest of the result follows from dynamic programming. The conclusion of the result follows from the fact that we have
  \begin{equation}
 V(0,X(0))= ||z||_{L_2([0,T])}^2-\gamma^2||w^*||_{L_2([0,T])}^2+\int_0^Tr(s,u^*(s))\ds+\E[q_T(X(T))],
\end{equation}
which is equivalent to say that
\begin{equation}
 ||z||_{L_2([0,T])}^2=\gamma^2||w^*||_{L_2([0,T])}^2-\int_t^Tr(s,u^*(s))\ds-\E[q_T(X(T))]+V(0,X(0))
\end{equation}
Under zero initial conditions and no terminal cost, this reduces to
\begin{equation}
 ||z||_{L_2([0,T])}^2=\gamma^2||w^*||_{L_2([0,T])}^2-\int_t^Tr(s,u^*(s))\ds<\gamma^2||w||_{L_2([0,T])}^2
\end{equation}
and the inequality follows from the fact that $r$ is nonnegative.
\end{proof}

We now address the important subclass of stochastic reaction networks with mass-action kinetics considered with a quadratic running cost:
\begin{hyp}\label{hyp:CTHinf3}
  The propensity functions of the network \net are of the form
  \begin{equation}\label{eq:propensities_bimolecular2}
  \lambda_i(t,x,u,w)=f_i(t,x)+g_{u,i}(t,x)u+g_{w,i}(t,x)w,\ i=1,\ldots,K
\end{equation}
where the $f_i(t,x)$'s are suitable polynomials and the $g_{u,i}(t,x)$'s and $g_{w,i}(t,x)$'s are suitable row-vector polynomials. Moreover, the function $r(t,u)$ in the cost \eqref{eq:costCTgeneral} takes the form
\begin{equation}\label{eq:kdslkd;skdkaldkal2}
    r(t,u)=u^{\T}R(t)u
  \end{equation}
  where $R(\cdot)\succ0$ and the controlled output is such that $z=h(t,x)$.
\end{hyp}

This leads to the following result:
\begin{theorem}
  Let us consider a stochastic reaction network satisfying Assumption \ref{hyp:CTHinf3}. Then an optimal pair $(u^*,w^*)$ that solves the optimization problem \eqref{eq:Hinfminimax} is given by
\begin{equation}\label{eq:Hinfoptimalpoint2}
\begin{array}{rcl}
  u^*(t)&=&-\dfrac{1}{2}R(t)^{-1}g_u(t,x)^{\T}\Delta V(t,x),\\
  w^*(t)&=&\dfrac{1}{2\gamma}g_w(t,x)^{\T}\Delta V(t,x),
\end{array}
\end{equation}
where the function $V(t,x)$ satisfies the Hamilton-Jacobi-Isaacs (HJI) equation given by
\begin{equation}\label{eq:HBICT2}
\begin{array}{l}
  V_t(t,x)+h(x)^{\T}h(x)+f(t,x)^{\T}\Delta V(t,x)-\dfrac{1}{4}\Delta V(t,x)^{\T}g_u(t,x)R(t)^{-1}g_u(t,x)^{\T}\Delta V(t,x)\\
  \qquad +\dfrac{1}{4\gamma^2}\Delta V(t,x)^{\T}g_u(t,x)g_u(t,x)^{\T}\Delta V(t,x)=0
\end{array}
\end{equation}
and $V(T,x)=q_T(x)$ where $\textstyle f(t,x)=\col_{i=1}^K(f_i(t,x))$, $\textstyle g_u(t,x)=\col_{i=1}^K(g_{u,i}(t,x))$, $\textstyle g_w(t,x)=\col_{i=1}^K(g_{w,i}(t,x))$, and $\textstyle \Delta V(t,x)=\col_{i=1}^K(\Delta V_i(t,x))$. Moreover, assuming zero initial conditions and no terminal cost, we have that
\begin{equation}
 ||z||_{L_2([0,T])}^2<\gamma^2||w||_{L_2([0,T])}^2
\end{equation}
for all $w\in \mathcal{W}_{0,T}$.
\end{theorem}
\begin{proof}
  A reaction network that satisfies Assumption \ref{hyp:CTHinf3} also satisfies Assumption \ref{hyp:CTHinf1} and Assumption \ref{hyp:CTHinf2}. The explicit expressions for \eqref{eq:Hinfoptimalpoint} can be obtained by differentiating \eqref{eq:HBICT} with respect to $u$ and $w$ and equating those expressions to zero.  Substitution of those explicit expressions in \eqref{eq:HBICT} yields \eqref{eq:HBICT2}.
\end{proof}}

\subsection{Control of unimolecular reaction networks}

\blue{In the unimolecular case, the propensity functions can be written as
\begin{equation}\label{eq:uni:prop:Hinf:CT}
\lambda_i(t,x,u,w)=\bar w_i(t)^{\T}\bar x+\bar b_i(t)^{\T}u+\bar{e}_i(t)^{\T}w,\ i=1,\ldots,K
\end{equation}
where $\bar w_i(t) = (w_{i,x}(t),w_{i,0}(t))$, $w_{i,x}\in\mathbb{R}_{\ge0}^{d}$, $\bar w_{i,0}(t)\in\mathbb{R}_{\ge0}$, $\bar b_i(t)\in\mathbb{R}_{\ge0}^{m}$, and $\bar{e}_i(t)\in\mathbb{R}_{\ge0}^{p}$, $i=1,\ldots,K$. The controlled output $z$ is assumed to depend linearly on the signals, that is
\begin{equation}\label{eq:uni:output:Hinf:CT}
  z(t)=C(t)X(t)+D(t)u(t)+F(t)w(t)+f(t)
\end{equation}
where $C:[0,T]\mapsto\mathbb{R}^{q\times d}$, $D:[0,T]\mapsto\mathbb{R}^{q\times m}$, $F:[0,T]\mapsto\mathbb{R}^{q\times p}$, and $f:[0,T]\mapsto\mathbb{R}^{q}_{\ge0}$ are assumed to be continuous functions. We also define $\bar C(t)=\begin{bmatrix}
  C(t) & f(t)
\end{bmatrix}$.

For this class of unimolecular networks, the moment equations read
\begin{equation}
\begin{array}{rcl}
\dfrac{\d}{\dt}\E[X(t)]&=&A(t)\E[X(t)]+B(t)\E[u(t)]+E(t)\E[w(t)]+d(t)\\
\E[z(t)]&=&C(t)\E[X(t)]+D(t)\E[u(t)]+F(t)\E[w(t)]+f(t)
\end{array}
\end{equation}
where
\begin{equation}
  A(t)=\sum_{k=1}^K\zeta_k \bar w_{k,x}(t)^{\T},\ B(t)=\sum_{k=1}^K\zeta_k \bar b_k(t)^{\T},\ E(t)=\sum_{k=1}^K\zeta_k\bar{e}_k(t)\ \textnormal{and }  d(t)=\sum_{k=1}^K\zeta_k \bar w_{k,0}(t)^{\T}.
\end{equation}

Finally, we consider the following quadratic version of the cost \eqref{eq:generalHinfcost} given by
\begin{equation}\label{eq:Hinf:cost:CT}
\begin{array}{rcl}
  J_T(u,w)&:=&\E\left[\displaystyle\int_0^T\left(u(s)^{\T}R(s)u(s)+||z(s)||_2^2-\gamma^2||w(s)||_2^2\right)\ds+\bar X(T)^{\T}Q_T\bar X(T)\right]
\end{array}
\end{equation}
where $\bar X=\col(X,1)$, $R:[0,T]\mapsto\pd^m$, $Q_T\in\psd^{d+1}$, and $\gamma>0$.
 We then have the following result:
\begin{theorem}\label{th:Hinf:uni}
  Assume that $D(\cdot)^{\T}F(\cdot)\equiv0$ and let $\gamma>0$. Let us further consider a unimolecular reaction network $(\X{},\mathcal{R})$ with propensity functions \eqref{eq:uni:prop:Hinf:CT} and controlled output \eqref{eq:uni:output:Hinf:CT}.  Then, the optimal pair $(u^*,w^*)$ that solves the optimization problem
  \begin{equation}
\min_{u\in\mathcal{U}_{0,T}}\max_{w\in\mathcal{W}_{0,T}}  J_T(u,w)
\end{equation}
  where $J_T(u,w)$ is given in \eqref{eq:Hinf:cost:CT} is given by
  \begin{equation}\label{eq:uopt:CT:Hinf}
  u^*(t,x(t))=-\dfrac{1}{2}(R(t)+D(t)^{\T}D(t))^{-1}\Omega_u(t)\bar x(t)
\end{equation}
and
\begin{equation}
  w^*(t,x(t))=(\gamma^2I-F(t)^{\T}F(t))^{-1}\Omega_w(t)\bar x(t)
\end{equation}
where 
\begin{equation}
  \begin{array}{rcl}
    \Omega_u(t)&:=&2D(t)^{\T}\bar C(t)+2\bar B(t)^{\T}P(t)+\sum_{k=1}^K\bar B_k(t)^{\T}P(t)\bar \zeta_ke_{d+1}^{\T}\\
    \Omega_w(t)&:=&2F(t)^{\T}\bar C(t)+2\bar E(t)^{\T}P(t)+\sum_{k=1}^K\bar E_k(t)^{\T}P(t)\bar \zeta_ke_{d+1}^{\T},
  \end{array}
\end{equation}
and where $P:[0,T]\mapsto\mathbb{S}^{d+1}$, $P(T)=Q_T$, solves the Riccati differential equation
\begin{equation}
\begin{array}{l}
  \dot{P}(t)+\bar A(t)^{\T}P(t)+P(t)\bar{A}(t)+\bar C(t)^{\T}\bar C(t)-\dfrac{1}{4}\Omega_u(t)^{\T}(R(t)+D(t)^{\T}D(t))^{-1}\Omega_u(t)\\
  \qquad+\dfrac{1}{4}\Omega_w(t)^{\T}(\gamma^2I-F(t)^{\T}F(t))^{-1}\Omega_w(t) +\dfrac{1}{2}\sum_{k=1}^K\He\left[\bar A_k(t)^{\T}P(t)\bar{\zeta}_ke_{d+1}^{\T}\right]=0
\end{array}
\end{equation}
where $\gamma^2I-F(t)^{\T}F(t)\succ0$, $\bar A_k(t)=\bar\zeta_k\bar w_k(t)^{\T}$, $\bar B_k(t)=\bar\zeta_k\bar b_k(t)^{\T}$, $\bar E_k(t)=\bar\zeta_k\bar e_k(t)^{\T}$,
\begin{equation}
  \bar A(t)=\begin{bmatrix}
    A(t) & d(t)\\
    0 & 0
  \end{bmatrix},\ \bar B(t)=\begin{bmatrix}
    B(t)\\
    0
  \end{bmatrix},\   \bar E(t)=\begin{bmatrix}
    E(t)\\
    0
  \end{bmatrix},\ \textnormal{and}\ \bar{C}(t)=\begin{bmatrix}
    C(t) & f(t)
  \end{bmatrix}.
\end{equation}
Moreover, in such a case, we have that
\begin{equation}
  ||z||_{L_2([0,T])}^2=\gamma^2||w^*||_{L_2([0,T])}^2-V(0,X(0))-||R^{1/2}u^*||_{L_2([0,T])}^2-||Q_T^{1/2}\bar{X}(T)||_2^2
\end{equation}
which reduces to
\begin{equation}
  ||z||_{L_2([0,T])}^2<\gamma^2||w||_{L_2([0,T])}^2
\end{equation}
for all $w\in\mathcal{W}_{0,T}$ whenever $X(0)=0$ and $Q_T=0$.
\end{theorem}

\begin{proof}
  In the present case, the Hamilton-Jacobi-Isaacs equation reads
  \begin{equation}\label{eq:djsjdlsadjsdlaj}
  V_t(t,x)+\min_u\max_w\left\{\sum_{i=1}^K\lambda_i(t,x,u,w)\Delta_iV(t,x)+u^{\T}R(t)u+z^{\T}z-\gamma^2w^{\T}w\right\}=0.
\end{equation}
First of all, note that
\begin{equation}
z(t)^{\T}z(t)=\begin{bmatrix}
    \bar x\\
    u\\
    w
  \end{bmatrix}^{\T}\mathcal{M}(t)\begin{bmatrix}
    \bar x\\
    u\\
    w
  \end{bmatrix}
\end{equation}
where
\begin{equation}
  \mathcal{M}(t):=\begin{bmatrix}
    \bar C(t)^{\T}\bar C(t) & \bar C(t)^{\T}D(t) & \bar C(t)^{\T}F(t)\\
    D(t)^{\T}\bar C(t) & D(t)^{\T}D(t) & 0\\
    F(t)^{\T}\bar C(t) & 0 & F(t)^{\T}F(t)
  \end{bmatrix}
\end{equation}
and where we have used the fact that $D(t)^{\T}F(t)=0$. Then, computing the derivative of the expression in the min-max expression in \eqref{eq:djsjdlsadjsdlaj} with respect to $u$ and $w$, and equating those expressions to zero yields
\begin{equation}
  \begin{array}{rcl}
    \sum_{i=1}^K\bar{b}_i^{\T}\Delta_iV(t,x)+2\bar x^{\T}\bar C(t)^{\T}D(t)+2u^{\T}(R(t)+D(t)^{\T}D(t))&=&0,\\
    \sum_{i=1}^K\bar{e}_i^{\T}\Delta_iV(t,x)+2\bar{x}^{\T}\bar C(t)^{\T}F(t)+2w^{\T}(F(t)^{\T}F(t)-\gamma^2I)&=&0.
  \end{array}
\end{equation}
Solving those equations for $u$ and $w$ yields
\begin{equation}
  \begin{array}{rcl}
    u^*&=&-\dfrac{1}{2}(R(t)+D(t)^{\T}D(t))^{-1}\left[2D(t)^{\T}\bar C(t)\bar x+\sum_{i=1}^K\bar{b}_i\Delta_iV(t,x)\right],\\
    w^*&=&\dfrac{1}{2}(\gamma^2I-F(t)^{\T}F(t))^{-1}\left[2F(t)^{\T}\bar C(t)\bar x+\sum_{i=1}^K\bar{e}_i\Delta_iV(t,x)\right]
  \end{array}
\end{equation}
where we have used the fact that $\gamma^2I-F(t)^{\T}F(t)\succ0$, and hence invertible, since $\gamma^2||w||_{L_{2,T}}^2\ge||z||_{L_{2,T}}^2>||Fw||_{L_{2,T}}^2$. Similarly, $R(t)+D(t)^{\T}D(t)$ is invertible since $R(t)+D(t)^{\T}D(t)\succeq R(t)\succ0$. The Hessian at $(u^*,w^*)$ is given by
\begin{equation}
  \begin{bmatrix}
    R(t)+D(t)^{\T}D(t) & 0\\0 & F(t)^{\T}F(t)-\gamma^2I
  \end{bmatrix}
\end{equation}
which proves that we have a minimum at $u=u^*$ since $R(t)+D(t)^{\T}D(t) \succ0$ and a maximum at $w=w^*$ since $F(t)^{\T}F(t)-\gamma^2I\prec0$.
 We also have that
\begin{equation}
\begin{array}{rcl}
  \sum_{i=1}^K\bar{b}_i\Delta_iV(t,x)&=&    \sum_{i=1}^K\left[2\bar{b}_i\bar\zeta_i^{\T}P(t)\bar x+b_i\bar\zeta_i^{\T}P(t)\bar \zeta_i\right]\\
                                                                     &=&    2\bar B(t)^{\T}P(t)\bar x+\sum_{i=1}^K\bar E_i^{\T}P(t)\bar \zeta_i\\
                                                                     &=&    \left(2\bar B(t)^{\T}P(t)+\sum_{i=1}^K\bar B_i^{\T}P(t)\bar \zeta_ie_{d+1}^{\T}\right)\bar x.
\end{array}
\end{equation}
Similarly,
\begin{equation}
\begin{array}{rcl}
  \sum_{i=1}^K\bar{e}_i\Delta_iV(t,x)&=&\left(2\bar E(t)^{\T}P(t)+\sum_{i=1}^K\bar E_i^{\T}P(t)\bar \zeta_ie_{d+1}^{\T}\right)\bar x.
\end{array}
\end{equation}
Therefore, we get that
\begin{equation}
  \begin{array}{rcl}
    u^*(t)&=&-\dfrac{1}{2}(R(t)+D(t)^{\T}D(t))^{-1}\left[2D(t)^{\T}\bar C(t)+2\bar B(t)^{\T}P(t)+\sum_{i=1}^K\bar B_i^{\T}P(t)\bar \zeta_ie_{d+1}^{\T}\right]\bar x\\
    &=&:-\dfrac{1}{2}Z_u(t)\bar x,\\
    w^*(t)&=&\dfrac{1}{2}(\gamma^2I-F(t)^{\T}F(t))^{-1}\left[2F(t)^{\T}\bar C(t)+2\bar E(t)^{\T}P(t)+\sum_{i=1}^K\bar E_i^{\T}P(t)\bar \zeta_ie_{d+1}^{\T}\right]\bar x\\
    &=:&\dfrac{1}{2}Z_w(t)\bar x.
  \end{array}
\end{equation}

Hence, the following expression holds at optimality
\begin{equation}
\begin{array}{rcl}
  z^*(t)^{\T}z^*(t)-\gamma^2 w^*(t)^{\T}w(t)+u^{*}(t)^{\T}R(t)u^*(t)&=&\bar x^{\T}\begin{bmatrix}
    I\\
    -Z_u(t)/2\\
    Z_w(t)/2
  \end{bmatrix}^{\T}\mathcal{M}(t)\begin{bmatrix}
    I\\
    -Z_u(t)/2\\
    Z_w(t)/2
  \end{bmatrix}\bar x.
\end{array}
\end{equation}

As previously, we have that
\begin{equation}
  \begin{array}{rcl}
  \bar w_k(t)^{\T}\bar x\left(2\bar\zeta_k^{\T}P\bar x+\bar\zeta_k^{\T}P\bar\zeta_k\right) &=&\bar x^{\T}\left(\He\left[P(t)\bar\zeta_k\bar w_k(t)^{\T}+\dfrac{1}{2}\bar w_k(t)\bar\zeta_k^{\T}P(t)\bar\zeta_ke_{d+1}^{\T}\right]\right)\bar x\\
    &=&\bar x^{\T}\left(\He\left[P(t)\bar A_k(t)+\dfrac{1}{2}\bar A_k(t)^{\T}P(t)\bar{\zeta}_ke_{d+1}^{\T}\right]\right)\bar x
  \end{array}
\end{equation}
and, hence,
\begin{equation}
  \sum_{k=1}^K  \bar w_k(t)^{\T}\bar x\left(2\bar\zeta_k^{\T}P\bar x+\bar\zeta_k^{\T}P\bar\zeta_k\right)    =\bar x^{\T}\left(\He\left[P(t)\bar A(t)+\dfrac{1}{2}\sum_{k=1}^K\bar A_k(t)^{\T}P(t)\bar{\zeta}_ke_{d+1}^{\T}\right]\right)\bar x.
\end{equation}

We also have
\begin{equation}
\begin{array}{rcl}
  \sum_{k=1}^K\bar b_k(t)^{\T}u^{*}(t)\Delta_kV(t,x)&=&    u^{*}(t)^{\T}\sum_{k=1}^K\bar b_k(t)\Delta_kV(t,x)\\
                                                                                  &=&   u^{*}(t)^{\T}\left(2\bar B(t)^{\T}P(t)+\sum_{i=1}^K\bar B_i(t)^{\T}P(t)\bar \zeta_ie_{d+1}^{\T}\right)\bar x\\
                                                                                  &=&  -\dfrac{1}{2}\bar x^{\T}Z_u(t)^{\T}\left(2\bar B(t)^{\T}P(t)+\sum_{i=1}^K\bar B_i(t)^{\T}P(t)\bar \zeta_ie_{d+1}^{\T}\right)\bar x
\end{array}
\end{equation}
Similarly,
\begin{equation}
\begin{array}{rcl}
  \sum_k\bar e_k(t)^{\T}w^*(t)\Delta_kV(t,x) &=&  \dfrac{1}{2}\bar x^{\T}Z_w(t)^{\T}\left(2\bar E(t)^{\T}P(t)+\sum_{i=1}^K\bar E_i(t)^{\T}P(t)\bar \zeta_ie_{d+1}^{\T}\right)\bar x
\end{array}
\end{equation}

Summing all those terms yields the following explicit expression for \eqref{eq:djsjdlsadjsdlaj}
\begin{equation}
\begin{array}{l}
  \bar x^{\T}\left[\dot{P}(t)+P(t)\bar A+\bar A^{\T}P(t)+\bar C(t)^{\T}\bar C(t)+\dfrac{1}{4}Z_u(t)^{\T}(D(t)^{\T}D(t)+R(t))Z_u(t)\right.\\
  \qquad +\dfrac{1}{4}Z_w(t)^{\T}(F(t)^{\T}F(t)-\gamma^2I)Z_w(t)-Z_u(t)^{\T}D(t)^{\T}\bar C(t)+Z_w(t)F(t)^{\T}\bar C(t) \\
  \qquad -\dfrac{1}{2}Z_u(t)^{\T}\left(2\bar B(t)^{\T}P(t)+\sum_{i=1}^K\bar B_i(t)^{\T}P(t)\bar \zeta_ie_{d+1}^{\T}\right)\\
  \qquad  \left.-\dfrac{1}{2}Z_w(t)^{\T}\left(2\bar E(t)^{\T}P(t)+\sum_{i=1}^K\bar E_i(t)^{\T}P(t)\bar \zeta_ie_{d+1}^{\T}\right)+\dfrac{1}{2}\sum_{k=1}^K\He[\bar A_k(t)^{\T}P(t)\bar{\zeta}_ke_{d+1}^{\T}]\right]\bar x=0
\end{array}
\end{equation}
which simplifies to
\begin{equation}
\begin{array}{l}
  \bar x^{\T}\left[\dot{P}(t)+P(t)\bar A+\bar A^{\T}P(t)+\bar C(t)^{\T}\bar C(t)-\dfrac{1}{2}Z_u(t)^{\T}(D(t)^{\T}D(t)+R(t))Z_u(t)\right.\\
  \\ \quad \left.+\dfrac{1}{4}Z_w(t)^{\T}(\gamma^2I-F(t)^{\T}F(t))Z_w(t)+\dfrac{1}{2}\sum_{k=1}^K\He[\bar A_k(t)^{\T}P(t)\bar{\zeta}_ke_{d+1}^{\T}]\right]\bar x=0.
\end{array}
\end{equation}
Since this must hold for all $\bar x=\col(x,1)$ with $x\in\mathbb{Z}_{\ge0}^d$, then this is equivalent to saying that
\begin{equation}
\begin{array}{l}
  \dot{P}(t)+P(t)\bar A(t)+\bar A(t)^{\T}P(t)+\bar C(t)^{\T}\bar C(t)-\dfrac{1}{2}Z_u(t)^{\T}(D(t)^{\T}D(t)+R(t))Z_u(t)\\
  \qquad +\dfrac{1}{4}Z_w(t)^{\T}(\gamma^2I-F(t)^{\T}F(t))Z_w(t)+\dfrac{1}{2}\sum_{k=1}^K\He[\bar A_k(t)^{\T}P(t)\bar{\zeta}_ke_{d+1}^{\T}]=0,
\end{array}
\end{equation}
which proves the result.
\end{proof}

As for the LQ control, the Riccati equation obtained for stochastic reaction networks is quite different that the deterministic counterpart. The exact closed-loop $H_\infty$-norm (or, equivalently, the $L_2$-gain) can be computed using a simple bisection algorithm which checks, at each iteration, whether a solution to the Riccati differential equation exists. Indeed, if the chosen $\gamma$ is smaller than the minimum achievable closed-loop $H_\infty$-norm, then the Riccati equation will have a finite escape time.}

\subsection{$H_\infty$ analysis of stochastic unimolecular reaction networks}

\blue{It seems interesting to address the $H_\infty$ analysis problem for completeness. To this aim, we consider the unimolecular stochastic reaction network \net with propensity functions
\begin{equation}\label{eq:uni:prop:Hinf:CT:OL}
\lambda_i(t,x,w)=\bar w_i(t)^{\T}\bar x+\bar{e}_i(t)^{\T}w,\ i=1,\ldots,K,
\end{equation}
and controlled output
\begin{equation}\label{eq:dksalkdas;dkl;kd;}
  z(t)=C(t)X(t)+F(t)w(t).
\end{equation}
This leads us to the following result:
\begin{theorem}[Bounded Real Lemma]
  Let $\gamma$ be a positive scalar. Let us further consider a unimolecular reaction network $(\X{},\mathcal{R})$ with propensity functions \eqref{eq:uni:prop:Hinf:CT:OL} and controlled output \eqref{eq:dksalkdas;dkl;kd;}. Then, this reaction network has finite-time $L_2$-gain if and only if there exists a matrix-valued function $P:[0,T]\mapsto\mathbb{S}^{d+1}$, $P(T)=0$, that solves the Riccati differential equation
\begin{equation}\label{eq:HinfOL}
\begin{array}{l}
  \dot{P}(t)+\bar A(t)^{\T}P(t)+P(t)\bar{A}(t)+\bar C(t)^{\T}\bar C(t)+\dfrac{1}{2}\sum_{k=1}^K\He\left[\bar A_k(t)^{\T}P(t)\bar{\zeta}_ke_{d+1}^{\T}\right]\\
  \qquad+\dfrac{1}{4}\Omega_w(t)^{\T}(\gamma^2I-F(t)^{\T}F(t))^{-1}\Omega_w(t)=0
\end{array}
\end{equation}
where $\gamma^2I-F(t)^{\T}F(t)\succ0$, $\bar A_k(t)=\bar\zeta_k\bar w_k(t)^{\T}$, $\bar E_k=\bar\zeta_k\bar e_k(t)^{\T}$, and
\begin{equation}
    \Omega_w(t):=2F(t)^{\T}\bar C(t)+2\bar E(t)^{\T}P(t)+\sum_{k=1}^K\bar E_k(t)^{\T}P(t)\bar \zeta_ke_{d+1}^{\T},
\end{equation}
where the matrices are the same as those in Theorem \ref{th:Hinf:uni}. Moreover, in such a case, the worst-case disturbance input is given by
  \begin{equation}
  w^*(t,x(t))=(\gamma^2I-F(t)^{\T}F(t))^{-1} \Omega_w(t)\bar{x}(t),
\end{equation}
and we have that
\begin{equation}
  ||z||_{L_2([0,T])}^2\le\gamma^2||w||_{L_2([0,T])}^2
\end{equation}
for all $w\in\mathcal{W}_{0,T}$ whenever $X(0)=0$ and equality holds with $w=w^*$.
\end{theorem}
\begin{proof}
The proof is based on that of Theorem \ref{th:Hinf:uni} where we have set $\bar{B}(\cdot)=\bar{B}_k(\cdot)=0$, $D(\cdot)=0$. In such a case, the optimal control input defined in \eqref{eq:uopt:CT:Hinf} collapses to zero whereas the differential Riccati equation reduces to \eqref{eq:HinfOL}. Finally, assuming no terminal cost (i.e. $Q_T=0$), and zero initial conditions yield the result.
\end{proof}

This result can be considered as the Bounded Real Lemma for unimolecular reaction networks. In fact, it can be reformulated, in perhaps a more familiar form, in terms of the Linear Matrix Inequalities $\gamma^2I-F(t)^{\T}F(t)\succ0$ and
\begin{equation}
  \begin{bmatrix}
     \dot{P}(t)+\He\left[P(t)\bar{A}(t)+\dfrac{1}{2}\sum_{k=1}^K\bar A_k(t)^{\T}P(t)\bar{\zeta}_ke_{d+1}^{\T}\right]+\bar C(t)^{\T}\bar C(t) & \dfrac{1}{2}\Omega_w(t)^{\T}\\
     \dfrac{1}{2}\Omega_w(t) & F(t)^{\T}F(t)-\gamma^2I
  \end{bmatrix}\preceq0
\end{equation}
which must hold for all $t\in[0,T]$. An approximate solution can then be found by using polynomial optimization methods such as SOS programming; see e.g. \cite{Parrilo:00,sostools3}.}

\section{Optimal Sampled-data $H_\infty$ Control of Stochastic Reaction Networks}\label{sec:SD:Hinf}

We now generalize the results of the previous sections to derive an $H_\infty$ control strategy in the sampled-data setting. Again, we exploit the hybrid formulation to represent the sampled-data controlled reaction network with uncertain norm-bounded inputs. Besides the already discussed difficulties in discretizing the dynamics of the stochastic reaction network in Section \ref{sec:CT:Hinf}, it is important to also add here that the hybrid formulation allows one to consider the input and output function spaces without needing to consider integral operators as required, for instance, in the lifting approach; see e.g. \cite{Chen:95}. The hybrid formulation allows one to deal with the original system directly and to rely on tools from hybrid systems theory. As previously, a general result is given together with its specialization to unimolecular networks.

\subsection{Preliminaries}

As previously, we consider that the network $(\X{},\mathcal{R})$ has $m$ control inputs and $p$ exogenous inputs. This motivates the following assumption

\begin{hyp}\label{hyp:SDHinf}
  The propensity functions of the network \net are of the form $ \lambda(t,x,v,w_c)$ where $v\in\mathbb{R}^m_{\ge0}$ is the state of the zero-order hold and $w_c\in\mathbb{R}_{\ge0}^{p_c}$ is the vector of continuous  exogenous disturbances. We also define the controlled outputs $z_c(t)=h_c(t,x(t),v(t),w_c(t))$ and $z_d(k)=h_d(t_k,x(t_k),v(t_k),w_d(k))$  where $w_d\in\mathbb{R}_{\ge0}^{p_d}$ is the vector of discrete exogenous disturbances and
\begin{equation}\label{eq:outputs:HinfSD}
  \begin{array}{rcl}
  h_c: \mathbb{R}_{\ge0}\times \mathbb{Z}_{\ge0}^d\times\mathbb{R}^m_{\ge0}\times\mathbb{R}_{\ge0}^{p_c}&\mapsto&\mathbb{R}^{q_c},\\
  h_d: \mathbb{R}_{\ge0}\times \mathbb{Z}_{\ge0}^d\times\mathbb{R}^m_{\ge0}\times\mathbb{R}_{\ge0}^{p_d}&\mapsto&\mathbb{R}^{q_d}
  \end{array}
  \end{equation}
  are assumed to be continuous functions.
\end{hyp}
%
%

We extend here the notion of $L_2$-norm to hybrid signals:
\begin{define}
The $L_2\times\ell_2$-norm of a hybrid signal $w:=(w_c,w_d)$ where $w_c:[0,T]\mapsto\mathbb{R}_{\ge0}^{p_c}$ and $w_d:\{0,\ldots,N\}\mapsto\mathbb{R}_{\ge0}^{p_d}$ is defined as
\begin{equation}
  \left|\left|\begin{bmatrix}
    w_c\\
    w_d
  \end{bmatrix}\right|\right|_{L_2([0,T])\times\ell_2([0,N])}:=\left(\int_0^T\E[||w_c(s)||_2^2]\ds+\sum_{k=0}^{N}\E[||w_d(k)||_2^2]\right)^{1/2}.
\end{equation}
\end{define}

We also need the following extension of the $L_2$-gain:
\begin{define}
Let us consider the reaction network \net which satisfies Assumption \ref{hyp:SDHinf} and assume that it defines an operator $L_2([0,T],\mathbb{R}_{\ge0}^{p_c})\times \ell_2(\{0,\ldots,N\},\mathbb{R}_{\ge0}^{p_d})\ni (w_c,w_d)=:z\mapsto w:=(z_c,z_d)\in L_2([0,T],\mathbb{R}_{\ge0}^{q_c})\times  \ell_2(\{0,\ldots,N\},\mathbb{R}_{\ge0}^{p_d})$. Then, the finite-horizon  $L_2\times\ell_2$-gain of that operator is defined as
\begin{equation}
    \left|\left|w\mapsto z\right|\right|_{L_2([0,T])\times\ell_2([0,N])-L_2([0,T])\times\ell_2([0,N])}:=\sup_{\substack{||w||_{L_2([0,T])\times\ell_2([0,N])}=1\\X(0)=0}}\left|\left|z\right|\right|_{L_2([0,T])\times\ell_2([0,N])}
\end{equation}
for all adapted signals $w$ in $L_2([0,T])\times\ell_2([0,N])$.
\end{define}

Based on the above definitions, we can state the main problem of the section which is to solve the following optimization problem
\begin{equation}
  \min_u   \left|\left|\begin{bmatrix}
    w_c\\
    w_d
  \end{bmatrix}\mapsto \begin{bmatrix}
    z_c\\
    z_d
  \end{bmatrix}\right|\right|_{L_2([0,T])\times\ell_2([0,N])-L_2([0,T])\times\ell_2([0,N])}
\end{equation}
where the control input $u$ is adapted and independent of $w$. In order to solve this problem, the following cost is introduced:
\begin{equation}\label{eq:HinfSDcost}
\begin{array}{rcl}
  J_T(u,w_c,w_d)&:=&\E\left[\displaystyle\int_0^T\left(||z_c(s)||_2^2-\gamma^2||w_c(s)||_2^2\right)\ds+q_T(\bar X(T),v(T))\right.\\
  &&\left.+\sum_{k=0}^{N-1}\left(r(t_k,u(t_k))+||z_d(k)||_2^2-\gamma^2||w_d(k)||_2^2\right)\right]
\end{array}
\end{equation}
where $r:[0,T]\times \mathbb{R}_{\ge0}^m\mapsto\mathbb{R}_{\ge0}$ is the running cost and $q_T:\mathbb{Z}_{\ge0}^d\mapsto\mathbb{R}$ is the terminal cost. The control input is chosen in such a way that it minimizes the worst case cost
\begin{equation}
  J_T^*=\min_{u\in\mathcal{V}_{0,T}} \max_{w_c\in\mathcal{W}^c_{0,T},w_d\in\mathcal{W}^d_{0,T}} J_T(u,w_c,w_d)
\end{equation}
where we have assumed that the minimum and maximum both exist, and
\begin{equation}
  \begin{array}{rcl}
   \mathcal{V}_{s,t}&:=&\left\{u\in \ell_2([s,t]\cap\{t_0,\ldots,t_N\}\times\mathbb{Z}_{\ge0}^d,\mathbb{R}_{\ge0}^m):\ u\ \textnormal{mesurable}\right\},\\
    \mathcal{W}^c_{s,t}&:=&\left\{w_c\in L_2([s,t]\times\mathbb{Z}_{\ge0}^d\times\mathbb{R}^m_{\ge0},\mathbb{R}_{\ge0}^{p_c}):\ w_c\ \textnormal{mesurable}\right\},\\
    \mathcal{W}^d_{s,t}&:=&\left\{w_d\in \ell_2([s,t]\cap\{t_0,\ldots,t_N\}\times\mathbb{Z}_{\ge0}^d\times\mathbb{R}^m_{\ge0},\mathbb{R}_{\ge0}^{p_d}):\ w_d\ \textnormal{mesurable}\right\},
  \end{array}
\end{equation}
for all $T\ge t\ge s\ge 0$. To solve this problem, we rely, once again, on Dynamic Programming and we introduce the value function
\begin{equation}
\begin{array}{rcl}
    V(t,x,v)& =&\min_{u\in\mathcal{V}_{t,T}} \max_{w_c\in\mathcal{W}^c_{t,T},w_d\in\mathcal{W}^d_{t,T}}\E\left[\int_t^T\left(||z(s)||_2^2-\gamma^2||w_c(s)||_2^2\right)\ds\right.\\
    &&\left.\left.+\sum_{k|t_k>t}^{N-1}\left(r(t_k,u(t_k))+||z_d(k)||_2^2-\gamma^2||w_d(k)||_2^2\right)\right|X(t)=x,v(t)=v\right]
\end{array}
\end{equation}
which satisfies $V(T,x,v)=q_T(x,v)$ for all $(x,v)\in\mathbb{Z}_{\ge0}^d\times\mathbb{R}_{\ge0}^m$. This yields the following result:
\begin{theorem}\label{th:kddsl;dls;dd;}
   Let us consider the sampled-data reaction network represented by \eqref{eq:SD:1}, \eqref{eq:SD:2}, \eqref{eq:SD:3}, which also satisfies Assumption \ref{hyp:SDHinf}. Then, the optimal control input and worst-case exogenous inputs are given by
   \begin{equation}
     \begin{array}{rcl}
       u^*(t_k)&=&\argmin_u \{V(t_k^+,x,u)+r(k,u)+z_d(k)^{\T}z_d(k)\},\\
       w_c^*(t)&=&\argmax_{w_c} \left\{\sum_{i=1}^K\lambda_i(t,x,v,w_c)\Delta V_i(t,x)+z_c(t)^{\T}z_c(t)-\gamma^2w_c(t)^{\T}w_c(t)\right\},\\
       w_d^*(k)&=&\argmax_{w_d} \{z_d(k)^{\T}z_d(k)-\gamma^2w_d(k)^{\T}w_d(k)\},
     \end{array}
   \end{equation}
   which hold for all $ t\in(t_k,t_{k+1}]$ and $k=0,\ldots,N-1$, and where the value function $V$ verifies the following  hybrid Hamilton-Jacobi-Isaacs (HJI) equation
\begin{equation}\label{eq:HBISD1}
  V_t(t,x,v)+\max_{w_c}\left\{\sum_{i=1}^K\lambda_i(t,x,v,w_c)\Delta_iV(t,x,v)+z_c^{\T}z_c-\gamma^2w_c^{\T}w_c\right\}=0,
\end{equation}
and
\begin{equation}\label{eq:HBISD2}
  \min_u\max_{w_d}\left\{V(t_k^+,x,u)-V(t_k,x,v)+r(t_k,u)+z_d^{\T}z_d-\gamma^2w_d^{\T}w_d\right\}=0,
\end{equation}
  which hold for all $ t\in(t_k,t_{k+1}]$ and $k=0,\ldots,N-1$, together with the terminal condition $V(T,x,v)=0$.
\end{theorem}

%
%

\subsection{Control of unimolecular reaction networks}

\blue{In this case of unimolecular networks, the propensity functions can be expressed as
\begin{equation}\label{eq:uni:prop:Hinf:SD}
  \lambda_i(t,x,v,w_c)=\bar{w}_i^{\T}\bar x+\bar{b}_i^{\T}v+\bar{e}_i^{\T}w_c,\ i=1,\ldots,K
\end{equation}
where $\bar{w}_i(\cdot)\in\mathbb{R}^{d+1}$, $\bar{b}_i(\cdot)\in\mathbb{R}^{m}$, and $\bar{e}_i(\cdot)\in\mathbb{R}^{p_c}$, $i=1,\ldots,K$. Similarly, we have that \eqref{eq:outputs:HinfSD} becomes
\begin{equation}
  \begin{array}{rcl}
    z_c(t)&=&C_c(t)X(t)+D_c(t)v(t)+F_c(t)w_c(t)+f_c(t),\\
    z_d(k)&=&C_d(k)X(t_k)+D_{d,v}(k)v(t_k)+D_{d,u}(k)u(t_k)+F_d(t_k)w_d(k)+f_d(k),
  \end{array}
\end{equation}
where the above matrix-valued functions are continuous and of appropriate dimensions.

Defining now the quadratic cost
\begin{equation}\label{eq:HinfSDcost}
\begin{array}{rcl}
  J_T(u,w_c,w_d)&:=&\E\left[\displaystyle\int_0^T\left(||z_c(s)||_2^2-\gamma^2||w_c(s)||_2^2\right)\ds+\tilde{x}(T)^{\T}Q_T\tilde{x}(T)\right.\\
  &&\left.+\sum_{k=0}^{N-1}\left(u(t_k)^{\T}R(k)u(t_k) +||z_d(k)||_2^2-\gamma^2||w_d(k)||_2^2\right)\right]
\end{array}
\end{equation}
allows us to formulate the following result:
\begin{theorem}
  Assume that $D_{d,u}(\cdot)^{\T}F_d(\cdot)\equiv0$ and let $\gamma>0$. Let us further consider a sampled-data unimolecular reaction network $(\X{},\mathcal{R})$ of the form \eqref{eq:SD:1}, \eqref{eq:SD:2} and \eqref{eq:SD:3} with propensity functions \eqref{eq:uni:prop:Hinf:SD}. Then, the optimal triplet $(u^*,w_c^*,w_d^*)$ that solves the optimization problem
\begin{equation}
  J_T^*=\min_{u\in\mathcal{V}_{0,T}} \max_{w_c\in\mathcal{W}^c_{0,T},w_d\in\mathcal{W}^d_{0,T}} J_T(u,w_c,w_d)
\end{equation}
where $J_T(u,w_c,w_d)$ is defined in \eqref{eq:HinfSDcost} is given by
\begin{equation}\label{eq:uopt:SD:Hinf}
  u^*(t_k,x,v)=-(P_3(t_k^+)+R(k)+D_{d,u}(k)^{\T}D_{d,u}(k))^{-1}\Omega_u(k)\begin{bmatrix}
    \bar x\\
    v
  \end{bmatrix}  
\end{equation}
\begin{equation}\label{eq:wcopt:SD:Hinf}
w_c^*(t,x,v)=(\gamma^2I-F_c(t)^{\T}F_c(t))^{-1}\Omega_c(t)\bar x
\end{equation}
and
\begin{equation}\label{eq:wdopt:SD:Hinf}
  w_d^*(k,x,v)=(\gamma^2I-F_d(k)^{\T}F_d(k))^{-1}F_d(k)^{\T}(C_d(k)\bar x+D_{d,v}(k)v)
\end{equation}
where 
\begin{equation}
  \begin{array}{rcl}
    \Omega_u(k)&:=&\begin{bmatrix}
      P_2(t_k^+)-C_d(k)^{\T}D_{d,u}(k)\\
      D_{d,v}(k)^{\T}D_{d,u}(k)
    \end{bmatrix}^{\T}\\
      \Omega_w(t)&:=&2F_c(t)^{\T}\bar C(t)+2\bar E(t)^{\T}P(t)+\sum_{k=1}^K\bar E_k(t)^{\T}P(t)\bar \zeta_ke_{d+1}^{\T},
  \end{array}
\end{equation}
where $P:[0,T]\mapsto\mathbb{S}^{d+m+1}$, $P(T)=Q_T$, solves the Riccati differential equation
\begin{equation}\label{eq:HinfSDCT}
\begin{array}{l}
  \dot{P}(t)+\bar A(t)^{\T}P(t)+P(t)\bar{A}(t)+\bar C(t)^{\T}\bar C(t) +\dfrac{1}{2}\sum_{k=1}^K\He\left[\bar A_k(t)^{\T}P(t)\bar{\zeta}_ke_{d+1}^{\T}\right]\\
  \qquad\qquad\qquad +\dfrac{1}{4}\Omega_w(t)^{\T}(\gamma^2I-F_c(t)^{\T}F_c(t))^{-1}\Omega_w(t)=0
\end{array}
\end{equation}
where $\gamma^2I-F_c(t)^{\T}F_c(t)\succ0$, $\bar A_k(t)=\bar\zeta_k\bar w_k(t)^{\T}$, $\bar B_k(t)=\bar\zeta_k\bar b_k(t)^{\T}$, $\bar E_k(t)=\bar\zeta_k\bar e_k(t)^{\T}$, and the Riccati difference equation
  \begin{equation}
    \begin{array}{l}
  \begin{bmatrix}
   I\\
   0
  \end{bmatrix} P_1(t_k^+)\begin{bmatrix}
   I\\
   0
  \end{bmatrix}^{\T}-P(t_k)+\bar{C}_d(k)^{\T}\bar{C}_d(k)+\bar{C}_d(k)^{\T}F_d(k)(\gamma^2I-F_d(k)^{\T}F_d(k))^{-1}F_d(k)^{\T}\bar{C}_d(k)\\
\qquad\qquad\qquad-\Omega_u(k)^{\T}(P_3(t_k^+)+R(k)+D_{d,u}(k)^{\T}D_{d,u}(k))^{-1}\Omega_u(k)=0
\end{array}
  \end{equation}
  where $\bar C_c(t):=\begin{bmatrix}
  C_c(t) & D_c(t) & f_c(t)
\end{bmatrix}$, $\bar C_d(k):=\begin{bmatrix}
  C_d(k) & D_{d,v}(k) & f_d(k)
\end{bmatrix}$,
  \begin{equation}
    \bar{A}(t):=\begin{bmatrix}
      A(t) & d(t) & \vline & B(t)\\
      0 & 0 & \vline & 0\\
      \hline
      0 & 0 & \vline & 0
    \end{bmatrix}, \bar{E}(t):=\begin{bmatrix}
      E(t)\\
      0\\
      \hline
      0
    \end{bmatrix}, \bar{\zeta}_k:=\begin{bmatrix}
      \zeta_k\\
      0\\
      \hline
      0
    \end{bmatrix},P(t)=\begin{bmatrix}
      P_1(t) & \vline & P_2(t)\\
      \hline
      P_2(t)^{\T} & \vline & P_3(t)
    \end{bmatrix},
  \end{equation}
  and $P_1(\cdot)\in\mathbb{S}^{d+1}$, $P_2(\cdot)\in\mathbb{R}^{(d+1)\times m}$, and $P_3(\cdot)\in\mathbb{S}^{m}$. Moreover, in such a case, we have that
\begin{equation}
\begin{array}{l}
   ||z_c||_{L_2([0,T])}^2+||z_d||_{\ell_2([0,N])}^2=\gamma^2(||w_c^*||_{L_2([0,T])}^2+||w_d^*||_{\ell_2([0,N])}^2)\\
   \qquad\qquad-V(0,X(0))-||R^{1/2}u^*||_{L_2([0,T])}^2-||Q_T^{1/2}\bar{X}(T)||_2^2
\end{array}
\end{equation}
which reduces to
\begin{equation}
||z_c||_{L_2([0,T])}^2+||z_d||_{\ell_2([0,N])}^2<\gamma^2(||w_c||_{L_2([0,T])}^2+||w_d||_{\ell_2([0,N])}^2)
\end{equation}
for all $w_c\in\mathcal{W}^c_{0,T}$, $w_c\in\mathcal{W}^d_{0,T}$, whenever $X(0)=0$ and $Q_T=0$.
  \end{theorem}
\begin{proof}
In the present case, the conditions of Theorem \ref{th:kddsl;dls;dd;} reduce to
\begin{equation}\label{eq:HBISD1:uni}
  V_t(t,x,v)+\max_w\left\{\sum_{i=1}^K\lambda_i(t,x,v,w)\Delta_iV(t,x,v)+z_c(t)^{\T}z_c(t)-\gamma^2w_c(t)^{\T}w_c(t)\right\}=0
\end{equation}
\begin{equation}\label{eq:HBISD2:uni}
  \min_u\max_w\left\{V_t(t_k^+,x,u)-V_t(t_k,x,v)+u^{\T}R(k)u+z_d(k)^{\T}z_d(k)-\gamma^2w_d(k)^{\T}w_d(k)\right\}=0
\end{equation}
and $V(T,x)=q_T(x)$. Computing the derivative with respect to $w_c$ in \eqref{eq:HBISD1:uni} yields
\begin{equation}
  \sum_{i=1}^K\bar{e}_i(t)^{\T}\Delta_iV(t,x,v)+2\bar x^{\T}C_c(t)^{\T}F_c(t)+2v^{\T}D_c(t)^{\T}F_c(t)+2w_c^{\T}F_c(t)^{\T}F_c(t)-2\gamma^2w_c^{\T}.
\end{equation}
Equating this to zero and solving for $w_c$ yields
\begin{equation}
  w_c^*(t)=\dfrac{1}{2}(\gamma^2I-F_c(t)^{\T}F_c(t))^{-1}\left[\sum_{k=1}^K\bar e_k(t)\Delta_kV(t,x,v)+2F_c(t)^{\T}\bar C_c(t)\tilde{x}(t)\right]
\end{equation}
where $\tilde x=(x,1,v)$. This is, indeed, a maximum if and only if $F_c^{\T}F_c-\gamma^2$ is negative definite. Substituting this expression in \eqref{eq:HBISD1:uni} yields the Riccati differential equation \eqref{eq:HinfSDCT}.

%

Analogously, the derivative of \eqref{eq:HBISD2:uni} with respect to $w_d$ is given by
\begin{equation}
  2\bar x^{\T}P_2(t_k^+)+2u^{\T}P_3(t_k^+)+2u^{\T}R(k)+2\bar x^{\T}C_d(k)^{\T}D_{d,u}(k)+2v^{\T}D_{d,v}(k)^{\T}D_{d,u}(k)+2u^{\T}D_{d,u}(k)^{\T}D_{d,u}(k)
\end{equation}
and equating it to 0 yields
\begin{equation}
\begin{array}{rcl}
  u^*(t_k,\bar  x,v)&=&-(P_3(t_k^+)+R(k)+D_{d,u}(k)^{\T}D_{d,u}(k))^{-1}(P_2(t_k^+)^{\T}\bar x+D_{d,u}^{\T}C_d(k)\bar x+D_{d,u}(k)^{\T}D_{d,v}(k)v)\\
  &=:&K_x(k)\bar{x}+K_v(k)v:=K(k)\tilde x.
\end{array}
\end{equation}
Similar calculations for $w_d$ yield
\begin{equation}
  2\bar x^{\T}C_d(k)^{\T}F_d(k)+2v^{\T}D_{d,v}(k)^{\T}F_d(k)+2w_d^{\T}F_d(k)^{\T}F_d(k)-2\gamma^2w_d^{\T}
\end{equation}
which leads to the optimal value
\begin{equation}
  w_d^*(k)=(\gamma^2I-F_d(k)^{\T}F_d(k))^{-1}F_d(k)^{\T}(C_d(k)\bar x+D_{d,v}(k)v)=W_x(k)\bar{x}+W_v(k)v:=W(k)\tilde x.
\end{equation}
The Hessian matrix is given by
\begin{equation}
  \begin{bmatrix}
    P_3(t_k^+)+R(k)+D_{d,u}(k)^{\T}D_{d,u}(k) & 0\\0 & F_d(k)^{\T}F_d(k)-\gamma^2
  \end{bmatrix}
\end{equation}
and shows that we get a maximum in $w_d$ if and only if $F_d^{\T}F_d-\gamma^2$ is negative definite where and a minimum in $u$ if and only if $P_3(t_k^+)+R(k)+D_{d,u}(k)^{\T}D_{d,u}(k)$ is positive definite. The latter condition is the case since $R(k)\succ0$ while the former is enforced by the condition in the result. Let $u^*=:K_x\bar x+K_vv$ and $w_d^*=:W_x\bar x+W_vv$. Substituting those expressions in \eqref{eq:HBISD2:uni} yields
\begin{equation}
\begin{array}{l}
  \begin{bmatrix}
    I & 0\\
    K_x(k) & K_v(k)
  \end{bmatrix}^{\T}P(t_k^+)\begin{bmatrix}
    I & 0\\
    K_x(k) & K_v(k)
  \end{bmatrix}-P(t_k)+K(k)^{\T}R(k)K(k)\\
  \qquad+\Omega_d(k)^{\T}\Omega_d(k)-\gamma^2W(k)^{\T}W(k)=0
\end{array}
\end{equation}
where $\Omega_d(k):=\bar{C}_d(k)+\bar{D}_{d,u}(k)K(k)+F_{d}(k)W(k)$. This can be reformulated as
\begin{equation}
\begin{array}{l}
  \begin{bmatrix}
    P_1(t_k^+) & 0\\
    0 & 0
  \end{bmatrix}+\He\left(\begin{bmatrix}
    P_2(t_k^+)\\
    0
  \end{bmatrix}K\right)+K(k)^{\T}(P_3(t_k^+)+R(k)+D_{d,u}^{\T}D_{d,u})K(k)\\
  \qquad+\He(\bar{C}_d^{\T}D_{d,u}K(k))+W(k)^{\T}(F_d(k)^{\T}F_d(k)-\gamma^2I)W(k)\\
  \qquad+\He(\bar{C}_d^{\T}F_{d}W(k))+\bar{C}_d^{\T}\bar{C}_d=0.
\end{array}
\end{equation}
Using now equalities
\begin{equation}
  \begin{array}{rcl}
    W(k)^{\T}(F_d(k)^{\T}F_d(k)-\gamma^2I)W(k)&=&-\bar{C}_d(k)^{\T}F_d(k)(\gamma^2I-F_d(k)^{\T}F_d(k))^{-1}F_d(k)^{\T}\bar{C}_d(k)\\
    \bar{C}_d(k)^{\T}F_dW(k)&=&\bar{C}_d(k)^{\T}F_d(k)(\gamma^2I-F_d(k)^{\T}F_d(k))^{-1}F_d(k)^{\T}\bar{C}_d(k).\\
    K(k)^{\T}(P_3(t_k^+)+R(k)+D_{d,u}^{\T}D_{d,u})K(k)&=&\Omega_d(k)^{\T}(P_3(t_k^+)+R(k)+D_{d,u}^{\T}D_{d,u})^{-1}\Omega_d(k)\\
    \left(\begin{bmatrix}
      P_2(t_k^+)\\
      0
    \end{bmatrix}+\bar{C}_d(k)^{\T}D_{d,u}(k)\right)K&=&-\Omega_d(k)^{\T}(P_3(t_k^+)+R(k)+D_{d,u}^{\T}D_{d,u})^{-1}\Omega_d(k)
  \end{array}
\end{equation}
where
\begin{equation}
  \Omega_d(k)=\begin{bmatrix}
      P_2(t_k^+)\\
      0
    \end{bmatrix}+\bar{C}_d(k)^{\T}D_{d,u}(k),
\end{equation}
yields the desired result.
\end{proof}}

\section{Concluding statements}

Optimal control results have been obtained for stochastic reaction networks. A difficulty of the approach is the impossibility of ensuring the nonnegativity of the control input, which results in the non-optimality of the control laws when restricted to nonnegative values. Such optimal control problems have been considered in the past  \cite{Heemels:98} but using the Pontryagin's Maximum Principle, which is not available in the current context. So, an alternative approach needs to be developed. Future extensions of the results would be the consideration of cost with infinite horizon through long-run average cost criterion.


\end{document}